\documentclass{amsart}
\usepackage{amssymb}
\usepackage{amsmath}
\usepackage{amsthm}
\usepackage{comment}
\usepackage{xcolor}
\usepackage{graphicx}
\usepackage{mathtools}
\usepackage{hyperref}
\usepackage{cleveref}
\usepackage{booktabs}
\usepackage{subfig}
\usepackage{float}
\usepackage{enumitem}

\newtheorem{theorem}{Theorem}[section]
\newtheorem{lemma}[theorem]{Lemma}
\newtheorem{corollary}[theorem]{Corollary}

\theoremstyle{definition}
\newtheorem{definition}[theorem]{Definition}

\theoremstyle{remark}
\newtheorem{remark}[theorem]{Remark}

\numberwithin{equation}{section}

\begin{document}

\title{A Nonlocal Biharmonic Model with $\Gamma$-Convergence to Local Model and an efficient numerical method}
\renewcommand*{\title}[2][]{\gdef\shorttitle{#1}\gdef\@title{#2}}
\title[A Nonlocal Biharmonic Model with $\Gamma$-Convergence]{}

\author{Weiye Gan}
\address{Department of Mathematical Sciences, Tsinghua University, Beijing, 100084, China.} \email{gwy22@mails.tsinghua.edu.cn}

\author{Tangjun Wang}
\address{Department of Mathematical Sciences, Tsinghua University, Beijing, 100084, China.} \email{wangtj20@mails.tsinghua.edu.cn}

\author{Qiang Du}
\address{Department of Applied Physics and Applied Mathematics, and Data Science Institute, Columbia
University, New York, NY 10027, USA.}  
\email{qd2125@columbia.edu}
\thanks{This work was supported by National Natural Science Foundation of China under grant 92370125 and US National Science Foundation DMS-2309245 and DMS-1937254.}

\author{Zuoqiang Shi}
\address{Yau Mathematical Sciences Center, Tsinghua University, Beijing, 100084, China \& Yanqi Lake Beijing Institute of Mathematical Sciences and Applications, Beijing, 101408, China.}
\email{zqshi@tsinghua.edu.cn}
\thanks{The forth author is the corresponding author.}

\subjclass[2020]{Primary 49J45, 65N30}

\keywords{nonlocal biharmonic model; Gamma convergence; nonlocal energy functional}

\begin{abstract}
  Nonlocal models and their associated theories have been extensively investigated in recent years. Among these, nonlocal versions of the classical Laplace operator have attracted the most attention, while higher-order nonlocal operators have been studied far less. In this work, we focus on the nonlocal counterpart of the classical biharmonic operator together with the clamped boundary condition ($u$ and $\frac{\partial u}{\partial n}$ are given on the boundary). We develop the variational formulation of a nonlocal biharmonic model, establish the existence and uniqueness of its solution, and analyze its convergence as the nonlocal horizon approaches zero. In addition, we propose an efficient finite element method to solve the nonlocal model and the numerical results verify the analytical properties of the nonlocal model and its solution.
\end{abstract}

\maketitle

\section{Introduction}
{Nonlocal models have found wide applications in different fields and have attracted much attention}~\cite{alfaro2017propagation,bavzant2002nonlocal,blandin2016well,dayal2007real,foss2022convergence,kao2010random, silling2000reformulation,vazquez2012nonlinear}. Among them, nonlocal diffusion models are of particular interest\cite{du2012analysis,du2019nonlocal,d2020numerical,du2013nonlocal,mengesha2013analysis,trask2019asymptotically,zhang2018accurate}
, as they {can be viewed as the continuum limit of discrete/graph Laplace operators and a} nonlocal version of the classical Laplace operator, both are widely used in many scientific and engineering disciplines. {There have been relatively limited studies of other nonlocal elliptic} operators, such as the nonlocal biharmonic operator.
{The latter will be the focus of our work here.}

{High-order elliptic differential operators have many applications. For example, in}  classical continuous mechanics, {the fourth-order biharmonic operator has been} used to model the deformation of thin plates, shells or beams. {When singularities, such as cracks, emerge, classical continuous models are no longer adequate. This motivates the need for multiscale modeling and the development of alternative continuum models such as  nonlcocal models}~\cite{silling2000reformulation,OGRADY20143177,OGRADY20144572}. Nonlocal fourth-order models are also used in semi-supervised learning, especially in cases with low label rate~\cite{cure2020}. When the sample rate is low, Laplace operator is not capable to properly distribute the label information while high-order operators give better results. 

{Based on the nonlocal formulation of the Laplace operator, it seems natural to derive a nonlocal biharmonic operator using two nonlocal Laplace operators of composition, which can be straightforward in the whole space or in periodic cells. This simple approach runs into complication in the case of bounded domains, due to the nonlocal nature of the operators.} 
For nonlocal diffusion model, many methods to impose Dirichlet or Neumann boundary conditions have been developed. {In particular, the notion of nonlocal volume constraints \cite{du2012analysis} is a natural and popular approach to replace the concept of local boundary conditions. 
Recently, many studies have been devoted to the construction of nonlocal models associated with different implementations of local or nonlocal conditions at or near the boundary}\cite{d2022prescription,yang2022uniform,lee2021second,shi2017convergence,shi2018harmonic,zhang2021second,scott2023nonlocal}.
In~\cite{gan2024convergence}, a framework is developed to approximate the local Dirichlet boundary condition using a nonlocal penalty formulation, while the exact Dirichlet boundary condition can be imposed by boundary localization \cite{scott2023nonlocal}. {For nonlocal biharmonic models with proper boundary conditions, studies remain limited.} In~\cite{tian2016class}, a general class of high-order nonlocal operators was proposed, the biharmonic operators were included as a special case. Based on the idea of volume constraint, nonlocal biharmonic operators with hinged and clamped boundary conditions were studied in~\cite{radu2017nonlocal}.

In this paper, we consider the nonlocal counterpart of the biharmonic equation with clamped boundary condition,
\begin{equation}
\label{eq:biharmonic}
\begin{aligned}
    &\Delta^2 u =f,\quad x\in \Omega\subset \mathbb{R}^d,\\
&u(x)=a(x),\quad \frac{\partial u}{\partial n}(x)=b(x),\quad x\in \partial \Omega.
\end{aligned}
\end{equation} {For the above local
biharmonic equation with the prescribd boundary conditions, the equivalent variational form is given by}
\begin{align*}
\min_{u\in H^2(\Omega)}\quad F(u)\,,
\end{align*}
{where the energy functional $F$ is defined} from $L^2(\Omega)$ to $ \mathbb{R}\cup \{+\infty\}$ as follows:
\begin{equation}\label{eq:local-functional}
			F(u)=\left\{\begin{array}{cc}
		\displaystyle	\int_\Omega|\Delta u|^2dx- 2 \int_\Omega fudx,&\text{if }u\in \mathcal{U},\\
				\infty&\text{otherwise}
			\end{array}\right.
		\end{equation}
where \begin{equation}\label{eq:lbh-space}
\mathcal{U}\coloneqq\left\{u\in H^2(\Omega)\big|u=a,\frac{\partial u}{\partial n}=b\text{ on }\partial\Omega\right\}.
	\end{equation}
for suitably specified boundary data $a$ and $b$ in appropriate trace spaces.
{As in \cite{tian2016class,radu2017nonlocal}, our idea is to use the corresponding nonlocal energy functional to  construct the resulting nonlocal model.} 
{For the Dirichlet energy term $\int_\Omega |\nabla u|^2dx$ in the energy $F$, one can use a standard formulation developed in the literature. Then, for the first term with the nonlocal Laplacian, we adopt a Green identity, similar to the study  point integral method~\cite{shi2017convergence}, to get the nonlocal approximation of the local Laplace operator together with the boundary condition on the normal derivative $\frac{\partial u}{\partial n}$.}
Finally, the boundary condition $u(x)=a(x),\; x\in \partial \Omega$ is enforced,  as in ~\cite{gan2024convergence}, by introducing a penalty term in the energy.
  {The nonlocal energy functional can then be derived based on the ideas presented above}. 

In addition, we propose an efficient numerical method for the nonlocal model based on finite element discretization. We consider the nonlocal model with Gaussian kernel and defined on tensor-product domain $\Omega=(0,1)^d$. In the finite element approach, multi-cubic polynomial basis is used. In the original finite element approach, to assemble the coefficient matrix, many $2d$ dimensional integrals need to be compute which impose huge computational load in the implementation. Taking advantage of the tensor-product structure of Gaussian kernel and multi-cubic polynomial basis, $2d$ dimensional integrals can be decomposed to be product of $2$ dimensional integrals which reduce the computational cost significantly. With the help of the efficient numerical solver, we conduct comprehensive study regarding the convergence of the nonlocal model as $\delta$
goes to 0. The numerical results shows that the nonlocal model has first order convergence to the counterpart local model. 

The remainder of the paper is organized as follows. {The specific form of the nonlocal energy functional and the main results of this paper are stated in Section 2, which include the well-posedness of the nonlocal models and their $\Gamma$-convergence as the nonlocal range $\delta\rightarrow 0$. The well-posedness is proved in Section 3.} In Section 4 and 5, $\Gamma$-convergence and convergence of the minimizers are proved, respectively. An efficient finite element method is presented in Section 6 and comprehensive numerical results are shown to validate the convergence of the nonlocal model. Some conclusions are drawn in Section 7.

	\section{Assumptions and main results}
	\label{sec:assump and result}
		Suppose that $\Omega$ is a bounded domain in $\mathbb{R}^d$ with $C^2$ boundary or a convex polytope, and $f\in L^2(\Omega)$. 
        $a$ and $b$ are the traces of a $H^2(\Omega)$ function and a normal derivative on $\partial \Omega$. $K$ and $R$ are nonnegative, monotonically decreasing, compactly supported one-dimensional $C^1$ kernels from $[0,\infty)$ to $[0,\infty)$. More specifically, we introduce the following assumptions on the kernels together with the notation for related constants.
  
  {\sc Assumption 2.1}\ \ Let $r_K$ and $r_R>0$ be given positive constants such that
		\begin{enumerate}
			\item[(K1)] $K,R:[0,\infty)\longrightarrow[0,\infty)$ belong to $C^1$;
			\item[(K2)] $K,R$ are monotonically decreasing;
			\item[(K3)] $ \mathrm{supp}(K)\subset[0,\frac{r_K^2}{4}
			]$, $\mathrm{supp}(R
			)\subset[0,\frac{r_R^2}{4}]$; 
   	\item[(K4)] $\displaystyle \bar{R}(x)\coloneqq \int_x^\infty R(s)ds, \; \forall\, x\geq 0$.
		\end{enumerate}

 For a constant $\delta>0$, we define the scaled kernel $R_\delta(x,y)\coloneqq \frac{1}{\delta^d}R(\frac{|x-y|^2}{4\delta^2})$, with $\bar{R}_\delta(x,y)$, $K_\delta(x,y)$ defined similarly.
		 
  The nonlocal model we construct is 
		\begin{equation}
			\label{eq:nonlocal-functional}
			\begin{aligned}
				F_n(u;a,b,f)=&\int_\Omega\left|\frac{1}{\delta_n^2}\int_{\Omega}R_{\delta_n}(x,y)(u(x)-u(y))dy-2\int_{\partial\Omega}\bar{R}_{\delta_n}(x,y)b(y)dS(y)\right|^2dx\\
				&-2\int_\Omega fudx+B_n(u,a)
			\end{aligned}
		\end{equation}
		for any $u\in L^2(\Omega)$ where 
		\begin{equation}
			\label{eq:boundary-term}
			B_n(u,a)\coloneqq\frac{1}{\xi_n}\int_{\partial\Omega}\left|\int_\Omega K_{\delta_n}(|x-y|)(a(x)-u(y))dy\right|^2dx
		\end{equation}
		is the boundary penalty term for Dirichlet boundary condition. $\xi_n$ is a positive constant satisfying $\lim_{n\rightarrow\infty}\xi_n=0$ and $\lim_{n\rightarrow\infty}\frac{\delta_n^2}{\xi_n}=0$, {see \cite{gan2024convergence}}. In the remainder of the paper, we denote $F_n(u)=F_n(u;a,b,f)$ to simplify the notation.
    \begin{remark}
        The nonlocal model enjoys a more concise form in the case of homogeneous boundary conditions. In terms of mathematical analysis, this could simply the presentation without sacrificing the generality. To extend such a case to the case of inhomogeneous boundary conditions, one may introduce the smooth extension $u_1$ of the boundary conditions, for example, by solving a biharmonic problem
        \begin{equation*}
            \begin{aligned}
    &\Delta^2 u =0,\quad x\in \Omega,\\
&u(x)=a(x),\quad \frac{\partial u}{\partial n}(x)=b(x),\quad x\in \partial \Omega.
\end{aligned}
        \end{equation*}
         With $u_1$ defined above, we could consider the function $u - u_1$ which satisfies homogeneous boundary conditions. However, such a construction might not be readily available in applications, and it might not be applicable in more general cases such as nonlinear problem. We thus choose to keep inhomogeneous boundary conditions in the discussion.
    \end{remark}
    \begin{remark}
        For the term $\int_{\Omega}fudx$ in \eqref{eq:nonlocal-functional}, we can use a mollified $f_{\delta_n}(x)\coloneqq\int_{\Omega}\tilde{R}_{\delta_n}(x,y)f(y)dy$ to replace $f(x)$. This change may relax the regularity constraints about $f$ and improve the regularity of the solution of \eqref{eq:nonlocal-functional} if $\tilde{R}$ is a carefully designed function of $R$, see \cite{wang2023nonlocal}
    \end{remark}
       
  Our main conclusions are the well-posedness of $F_n$ and the $\Gamma$-convergence from $F_n$ to $F$ as stated follows,
  
  \begin{theorem}[Well-posedness]
      \label{thm:well-posedness}
      Suppose that $\Omega$ is a bounded Lipschitz domain in $\mathbb{R}^d$. $K,R,\tilde{R}$ are kernels satisfying Assumption 2.1. Then, for arbitrary fixed $\delta_n$ and $\xi_n$, there exists a unique minimizer of $F_n:L^2(\Omega)\rightarrow\mathbb{R}$ defined as \eqref{eq:nonlocal-functional}.
  \end{theorem}
		\begin{theorem}[$\Gamma$-convergence]
			\label{thm:gamma-convergence}
			Suppose that $\Omega$ is a bounded domain in $\mathbb{R}^d$ with $C^2$ boundary or a convex polytope.  $K,R,\tilde{R}$ are kernels satisfying {{Assumption} 2.1}. $\{\delta_n\}$, $\{\xi_n\}$ are sequences of positive constants tending to 0 as $n\rightarrow\infty$ and satisfying $\lim_{n\rightarrow\infty}\frac{\delta_n^2}{\xi_n}=0$. Then we have
			\begin{equation*}
				F_n\stackrel{\varGamma}{\longrightarrow}F\quad\text{in }L^2(\Omega),
			\end{equation*}
			where $F_n,F$ are defined as \eqref{eq:nonlocal-functional},\eqref{eq:local-functional}.
		\end{theorem}

    \begin{theorem}[Convergence of minimizers]
	\label{thm:convergence-minimizers}
	Suppose that $\Omega$ is a bounded domain in $\mathbb{R}^d$ with $C^2$ boundary or a convex polytope.  $K,R,\tilde{R}$ are kernels satisfying {{Assumption} 2.1}. $\{\delta_n\}$, $\{\xi_n\}$ are sequences of positive constants tending to 0 as $n\rightarrow\infty$ and satisfying $\lim_{n\rightarrow\infty}\frac{\delta_n^2}{\xi_n}=0$. Then any sequence $\{u_n\}\subset L^2(\Omega)$ satisfying
        \begin{equation*}
            \lim_{n\rightarrow\infty}(F_n(u_n)-\inf_{u\in L^2(\Omega)}F_n(u))=0
	\end{equation*}
    is relatively compact in $L^2(\Omega)$ and
	\begin{equation*}
            \lim_{n\rightarrow\infty}F_n(u_n)=\min_{u\in L^2(\Omega)}F(u).
	\end{equation*}
	Furthermore, every cluster point of $\{u_n\}$ is a minimizer of $F$.
    \end{theorem}
    
  \section{Well-posedness of nonlocal functionals} In this section, we clarify the well-posedness of the minimization problem with respect to nonlocal functional $F_n$ defined as \cref{eq:nonlocal-functional}. 
  For technical reasons, we first define a nonlocal Dirichlet energy
  \begin{equation}	\label{eq:nonlocol-laplacian}
			L_n(u)\coloneqq\frac{1}{\delta_n^2}\int_\Omega\int_\Omega R_{\delta_n}(x,y)|u(x)-u(y)|^2dxdy.\\
		\end{equation}
    And the following lemma 
    is vital.
    \begin{lemma}
        \label{lem:poincare}
        Suppose that $\Omega$ is a bounded Lipschitz domain in $\mathbb{R}^d$. $\tilde{R},K,\{\delta_n\},$\\$\{\xi_n\}$ satisfy the same assumption as \cref{thm:gamma-convergence}. Then, there exists constants $C_1,C_2>0$ not depending on $n$ and $u$ such that
        \begin{equation*}
            L_n(u)+B_n(u,a)\geq C_1\left\lVert u\right\rVert_{L^2(\Omega)}^2-C_2\left\lVert a\right\rVert^{2}_{L^2(\partial\Omega)}
        \end{equation*}
        for all $u\in L^2(\Omega)$, where $L_n,B_n$ are defined as \eqref{eq:nonlocol-laplacian}, \eqref{eq:boundary-term}.
    \end{lemma}
    \begin{proof}
        Note that $\{\xi_n\}$ is bounded
        \begin{equation*}
            \begin{aligned}
&\quad B_n(u,a)\\
&\geq\frac{1}{\max_n\xi_n}\int_{\partial\Omega}\left|\int_\Omega K_{\delta_n}(|x-y|)(a(x)-u(y))dy\right|^2dx\\
                &\geq\frac{1}{\max_n\xi_n}\left\{\frac{1}{2}\int_{\partial\Omega}\left|\int_\Omega K_{\delta_n}(|x-y|)u(y)dy\right|^2dx-\int_{\partial\Omega}\left|\int_\Omega K_{\delta_n}(|x-y|)a(x)dy\right|^2dx\right\}\\
                &\geq \frac{1}{2\max_n\xi_n}\int_{\partial\Omega}\left|\int_\Omega K_{\delta_n}(|x-y|)u(y)dy\right|^2dx-C_2\left\lVert a\right\rVert^{2}_{L^2(\partial\Omega)}.
            \end{aligned}
        \end{equation*}
        Define 
        \begin{equation*}
            \tilde{u}(x)=\frac{1}{\int_\Omega K_{\delta_n}(|x-y|)dy}\int_\Omega K_{\delta_n}(|x-y|)u(y)dy.
        \end{equation*}
        With lemma 4.3, lemma 5.4 in~\cite{gan2024convergence} and the Poinc\'are inequality, we have 
        \begin{equation*}
        \begin{aligned}
            &\quad L_n(u)+\frac{1}{2\max_n\xi_n}\int_{\partial\Omega}\left|\int_\Omega K_{\delta_n}(|x-y|)u(y)dy\right|^2dx\\
            &\geq C\left\lVert\nabla \tilde{u}\right\rVert_{L^2(\Omega)}^2+\frac{1}{2\max_n\xi_n}\left\lVert \tilde{u}\right\rVert_{L^2(\partial\Omega)}^2\\
            &\geq C_0 \left\lVert \tilde{u}\right\rVert_{L^2(\Omega)}^2.
        \end{aligned}
        \end{equation*}
        Moreover, with the H\"older inequality, Lemma 4.3 in~\cite{gan2024convergence} and the fact that $\int_{\Omega}K_{\delta_n}(|x-y|)dy$ has a positive lower bound not depending on $n$, there exists a constant $C_3=C_3(K,R,\Omega)$, such that
        \begin{equation*}
        \begin{aligned}
            \left\lVert \tilde{u} - u\right\rVert_{L^2(\Omega)}^2\leq C_3L_n(u).
        \end{aligned}
        \end{equation*}
        Let $C_1 = C_0/(1+C_0C_3)$. we obtain the estimate.
    \end{proof}
    
    We also show that $L_n$ is controlled by the interior term and the boundary term in our nonlocal model $F_n$ 
   \begin{lemma}
       \label{lem:basic-control}Suppose that $\Omega$ is a bounded Lipschitz domain in $\mathbb{R}^d$. $\tilde{R},K,\{\delta_n\},$\\$\{\xi_n\}$ satisfy the same assumption as \cref{thm:gamma-convergence}. Then, there exists a constant $C>0$ not depending on $n$ and $u$ such that
        \begin{equation*}
            L_n(u)\leq C(F_n(u) + \int_{\Omega} fu dx+ \lVert a\rVert_{L^2(\partial\Omega)}^2+\lVert b\rVert_{L^2(\partial\Omega)}^2) 
        \end{equation*}
        for all $u\in L^2(\Omega)$, where $F_n,L_n$ are defined as \eqref{eq:nonlocal-functional}, \eqref{eq:nonlocol-laplacian}.
   \end{lemma}
   \begin{proof}
       Define 
       \begin{equation*}
           I_{n}(u) = \frac{1}{\delta_n^2}\int_{\Omega}R_{\delta_n}(x,y)(u(x)-u(y))dy-2\int_{\partial\Omega}\bar{R}_{\delta_n}(x,y)b(y)dS(y).
       \end{equation*}
       Then, 
       \begin{equation*}
           \begin{aligned}
               L_n(u)&=2\Big\langle u,\frac{1}{\delta_n^2}\int_\Omega R_{\delta_n}(x,y)(u(x)-u(y))dy\Big\rangle_{L_2(\Omega)}\\
               &=2\langle u, I_n(u)\rangle_{L_2(\Omega)} + 4\Big\langle u,\int_{\partial\Omega}\bar{R}_{\delta_n}(x,y)b(y)dS(y)\Big\rangle_{L_2(\Omega)}\\
               &\leq 2\lVert u\rVert_{L^2(\Omega)}\lVert I_n(u)\rVert_{L_2(\Omega)}+ 2\int_{\Omega}\int_{\partial\Omega}\bar{R}_{\delta_n}(x,y)u(x)^2dS(y)dx\\
               &+ 2\int_{\Omega}\int_{\partial\Omega}\bar{R}_{\delta_n}(x,y)b(y)^2dS(y)dx\\
               &\leq 2\lVert u\rVert_{L^2(\Omega)}\lVert I_n(u)\rVert_{L_2(\Omega)} + \xi_n C_1(K,R,\Omega)B_n(u,a)\\
               & + C_2(R,\Omega)(\lVert a\rVert_{L^2(\partial\Omega)}^2+\lVert b\rVert_{L^2(\partial\Omega)}^2)
           \end{aligned}
       \end{equation*}
       where the last inequality holds by applying Lemma 4.3 in \cite{gan2024convergence}. Then, note that
        \begin{equation*}
            F_n(u) + 2\int_\Omega fudx = \lVert I_n(u)\rVert_{L^2(\Omega)}^2 + B_n(u,a).
        \end{equation*}
        We obtain the conclusion with \cref{lem:poincare}.
   \end{proof}

    Note that $F_n$ is convex. We can demonstrate the existence of minimizers with a standard method for a convex functional. For convenience, we first give a definition that is one of the sufficient conditions. In this article, we only need to focus on Hilbert spaces rather than general Banach spaces.
    \begin{definition}[Ellipticity]
        Let $X$ be a Hilbert space. We call a functional $F:X\rightarrow\mathbb{R}\cup\{+\infty\}$ elliptic if $\lim_{\langle u,u\rangle_{X}\rightarrow\infty} F_n(u)=+\infty$
    \end{definition}
    \begin{lemma}[Corollary 3.23 in \cite{brezis2011functional}]
    \label{lem:existence-of-minimizers}
    Let $X$ be a Hilbert space. Suppose that a functional $F:X\rightarrow\mathbb{R}\cup\{+\infty\}$ is convex, continuous, and elliptic. Then, there exists some $x_0\in X$ such that
    \begin{equation*}
        F(x_0)=\min_X F.
    \end{equation*}
    \end{lemma}
    
    \begin{lemma}
    \label{lem:well-posedness1}
        Under the same assumption as \cref{thm:well-posedness}, $F_n:L^2(\Omega)\rightarrow\mathbb{R}$ defined as \cref{eq:nonlocal-functional} is convex, continuous, and elliptic.
    \end{lemma}
    \begin{proof}
        $F_n$ is convex since $f(x)=|x|^2$ is a convex function. With the basic inequality $\left\lvert\alpha^2-\beta^2\right\rvert\leq 2\left\lvert\alpha\right\rvert\left\lvert\alpha-\beta\right\rvert+(\alpha-\beta)^2$ and the H\"older inequality, one can demonstrate the continuity. For example, for fixed $u$ and arbitrary $u_1$, define $u_d = u_1-u$. Then,as $\left\lVert u_d\right\rVert_{L^2(\Omega)}\rightarrow 0$ we have
        \begin{equation*}
            \begin{aligned}
                \quad&\Bigg|\int_\Omega\left|\frac{1}{\delta_n^2}\int_{\Omega}R_{\delta_n}(x,y)(u(x)-u(y))dy-2\int_{\partial\Omega}\bar{R}_{\delta_n}(x,y)b(y)dS(y)\right|^2dx\\
                &-\int_\Omega\left|\frac{1}{\delta_n^2}\int_{\Omega}R_{\delta_n}(x,y)(u_1(x)-u_1(y))dy-2\int_{\partial\Omega}\bar{R}_{\delta_n}(x,y)b(y)dS(y)\right|^2dx\Bigg|\\
                &\leq C_1(\delta_n,R,\Omega)\int_\Omega\int_{\Omega}\left|u_d(x)-u_d(y)\right|^2dydx\\
                &+C_2(\delta_n,R,\bar{R},\Omega,b)\int_\Omega\left(\int_{\Omega}\left|u(x)-u(y)\right|dy+1\right)\int_{\Omega}\left|u_d(x)-u_d(y)\right|dydx\\
                &\leq C_3(\delta_n,R,\bar{R},\Omega,b)\left\lVert u_d\right\rVert_{L^2(\Omega)}\left(\left\lVert u\right\rVert_{L^2(\Omega)}+1\right)\longrightarrow 0.
            \end{aligned}
        \end{equation*}
        
        Moreover, for the ellipticity, with \cref{lem:poincare} and \cref{lem:basic-control}, there exists positive constants $C_4$, $C_5$ not depending on $u$ such that
        \begin{equation*}
            F_n(u)\geq C_4\left\lVert u\right\rVert_{L^2(\Omega)}^2-\left\lVert f\right\rVert_{L^2(\Omega)}\left\lVert u\right\rVert_{L^2(\Omega)}-C_5\longrightarrow+\infty
        \end{equation*}
        as $\left\lVert u\right\rVert_{L^2(\Omega)}\rightarrow\infty$.
    \end{proof}
    
    The Uniqueness is basically a direct consequence of the strong convexity of function $f(x)=x^2$ and coervicity established in \cref{lem:poincare}. Nevertheless, we still provide the proof for completeness and clarity.
    \begin{lemma}
        \label{lem:well-posedness2}
        Under the same assumption as \cref{thm:well-posedness}, suppose that $u_n^*$ is a minimizer of $F_n$ defined as \eqref{eq:nonlocal-functional}. Then, there exists a positive constant $C$ not depending on $u_n^*$ and $\delta_n$ such that
        \begin{equation*}
            F_n(u)-F_n(u_n^*)\geq C\left\lVert u-u_{n}^*\right\rVert_{L^2(\Omega)}^2
        \end{equation*}
        for all $u\in L^2(\Omega)$.
    \end{lemma}
    \begin{proof}
    Since 
    \begin{equation*}
        (\alpha a +(1-\alpha)b)^2=\alpha a^2 + (1-\alpha) b^2 -\alpha(1-\alpha)(a-b)^2
    \end{equation*}
    for all $a,b\in\mathbb{R}$ and $\alpha$ satisfying $0 \leq\alpha\leq 1.$ For any $u\in L^2(\Omega)$, define $v_n = u - u_n^*$. We have
    \begin{equation*}
    \begin{aligned}
        &\quad F_n(u_n^*)\leq F_n(\alpha u+(1-\alpha)u_n^*)\\
        &\leq \alpha F_n(u)+(1-\alpha)F_n(u_n^*)\\
        &-\alpha(1-\alpha)\left\{\int_{\Omega}\left|\frac{1}{\delta_n^2}\int_\Omega R_{\delta_n}(|x-y|)(v_n(x)-v_n(y))dy\right|^2dx + B_n(v_n,0)\right\}.
    \end{aligned}
    \end{equation*}
    Let $\alpha\rightarrow0^+$,
    \begin{equation*}
        F_n(u)-F_n(u_n^*)\geq\int_{\Omega}\left|\frac{1}{\delta_n^2}\int_\Omega R_{\delta_n}(|x-y|)(v_n(x)-v_n(y))dy\right|^2dx+ B_n(v_n,0).
    \end{equation*}
    With \cref{lem:poincare} and \cref{lem:basic-control} for $a,b\equiv 0$, there exists a constant $C$ not depending on $u_n^*$ and $\delta_n$ such that
    \begin{equation*}
        F_n(u)-F_n(u_n^*)\geq C\left\lVert u-u_{n}^*\right\rVert_{L^2(\Omega)}^2.
    \end{equation*}
    \end{proof}
    
    \cref{thm:well-posedness} is a direct conclusion with \cref{lem:existence-of-minimizers}, \cref{lem:well-posedness1} and \cref{lem:well-posedness2}. Moreover, when the boundary conditions are homogeneous, we can render the uniform boundedness for minimzers of $F_n$ with respect to different parameters $\delta_n,\xi_n$.

    \begin{corollary}
    \label{lem:boundedness-minimizers}
        Consider the homogeneous boundary conditions $a\equiv 0$ and $b\equiv 0$. Under the same assumption as \cref{thm:well-posedness}, suppose that $u_n^*$ is the minimizer of $F_n:L^2(\Omega)\rightarrow\mathbb{R}$ defined as \cref{eq:nonlocal-functional}. Then, there exists a constant $C>0$ not depending on $n$ such that
        \begin{equation*}
            \left\lVert u_n^*\right\rVert_{L^2(\Omega)}\leq C.
        \end{equation*}
    \end{corollary}
    \begin{proof}
    With calculus of variations, we can obtain the Euler-Lagrange equation of nonlocal model \eqref{eq:nonlocal-functional} with $a,b\equiv 0$ as follows,
    \begin{equation*}
    \begin{aligned}
		&\quad \frac{2}{\delta_n^4}\int_\Omega \int_\Omega R_{\delta_n}(x,y) \Big( R_{\delta_n}(x,z)(u_n^*(x)-u_n^*(z))-R_{\delta_n}(y,z)(u_n^*(y)-u_n^*(z))\Big) dz dy \\
		& +\frac{4\varepsilon}{\delta_n^2} \int_\Omega \tilde{R}_{\delta_n}(x,y) \left(u_n^*(x)-u_n^*(y)\right) dy\\
  &+ \frac2{\xi_n} \int_\Omega \left( \int_{\partial\Omega} K_{\delta_n}(x,z)K_{\delta_n}(y,z)dS(z) \right) u_n^*(y) dy   = f(x).
    \end{aligned}
    \end{equation*}
    Multiplying both sides by $u_n^*(x)$ and integrating on $\Omega$ with respect to $x$, we have
    \begin{equation*}
    \begin{aligned}
        &\quad\frac{2}{\delta_n^4}\int_\Omega\left(\int_\Omega R_{\delta_n}(x,y)\left(u_n^*(x)-u_n^*(y)\right)dy\right)^2dx\\
        &+\frac{4\varepsilon}{\delta_n^2} \int_\Omega\int_\Omega \tilde{R}_{\delta_n}(x,y) \left(u_n^*(x)-u_n^*(y)\right)^2 dydx+\frac2{\xi_n} \int_{\partial\Omega} \left(  \int_\Omega K_{\delta_n}(x,y)u_n^*(y) dy\right)^2dS(x)\\
        &=\int_{\Omega}f(x)u_n^*(x)dx.
    \end{aligned}
    \end{equation*}
    With \cref{lem:poincare}, there exists a constant $C_1>0$ not depending on $n$ such that the left side is not less than $C_1\left\lVert u_n^*\right\rVert_{L^2(\Omega)}^2$. With the H\"older inequality, we can select $C=\frac{\left\lVert f\right\rVert_{L^2(\Omega)}}{C_1}$ to satisfy the conclusion.
    \end{proof}

		\section{$\Gamma$-Convergence of nonlocal functionals}
		\label{sec:proof}
		In this section, we give the proof details of \cref{thm:gamma-convergence}. As mentioned in the introduction, the following estimate presented in \cite{shi2017convergence} is significant. The convergence rate of the residue $r_n$ is analyzed in \cite{shi2017convergence} for $H^3$ function. While we need to process lower regularity.
		\begin{lemma}
			\label{lem:appro}
			Suppose that $\Omega$ is a bounded Lipschitz domain in $\mathbb{R}^d$. $\{\delta_n\}$ is a sequence of positive constants tending to 0 as $n\rightarrow\infty$. Then, for all fixed $u\in H^2(\Omega)$, the $L^2(\Omega)$ norm of the following function
			\begin{equation*}
				\begin{aligned}
					r_n(x)&=\int_\Omega \bar{R}_{\delta_n}(x,y)\Delta u(y)dy+\frac{1}{\delta_n^2}\int_{\Omega}R_{\delta_n}(x,y)(u(x)-u(y))dy\\&-2\int_{\partial\Omega}\bar{R}_{\delta_n}(x,y)\frac{\partial u}{\partial n}(y)dS(y)
				\end{aligned}
			\end{equation*}
			tends to $0$ as $n\rightarrow\infty$.
   \end{lemma}
   
			\begin{proof}
				With Green formula,
				\begin{equation*}
					\begin{aligned}
						&\quad \int_{\partial\Omega}\bar{R}_{\delta_n}(x,y)\frac{\partial u}{\partial n}(y)dS(y)\\
						&=\int_\Omega \bar{R}_{\delta_n}(x,y)\Delta u(y)dy+\int_\Omega\nabla_y\bar{R}_{\delta_n}(x,y)\cdot\nabla u(y)dy\\
						&=\int_\Omega \bar{R}_{\delta_n}(x,y)\Delta u(y)dy-\frac{1}{2\delta_n^2}\int_\Omega R_{\delta_n}(x,y)\langle\nabla u(y),x-y\rangle dy.
					\end{aligned}
				\end{equation*}
				Denote by $\tilde{u}$ the extension of $u$ on $\mathbb{R}^d$. Then,
				\begin{equation*}
					\begin{aligned}
						r_n(x)&=-\int_\Omega \bar{R}_{\delta_n}(x,y)\Delta \tilde{u}(y)dy+\frac{1}{\delta_n^2}\int_\Omega R_{\delta_n}(x,y)(\tilde{u}(x)-\tilde{u}(y)+\langle\nabla \tilde{u}(y),x-y\rangle)dy\\
						&=r_1^n(x)+r_2^n(x)+r_3^n(x)
					\end{aligned}
				\end{equation*}
				where
				\begin{equation*}
					\begin{aligned}
						r_1^n(x)&\coloneqq\frac{1}{\delta_n^2}\int_\Omega R_{\delta_n}(x,y)\left(\tilde{u}(x)-\tilde{u}(y)+\langle\nabla \tilde{u}(y),x-y\rangle\right)dy\\
						&-\sum_{i,j=1}^d\frac{1}{\delta_n^2}\int_\Omega R_{\delta_n}(x,y)(x_i-y_i)(x_j-y_j)\int_0^1(1-t)\tilde{u}_{ij}(y+t(x-y))dtdy,
					\end{aligned}
				\end{equation*}
				\begin{equation*}
					\begin{aligned}
						r_2^n(x)\coloneqq\sum_{i\neq j}\frac{1}{\delta_n^2}\int_\Omega R_{\delta_n}(x,y)(x_i-y_i)(x_j-y_j)\int_0^1(1-t)\tilde{u}_{ij}(y+t(x-y))dtdy
					\end{aligned}
				\end{equation*}
				and
				\begin{equation*}
					\begin{aligned}
				r_3^n(x)&\coloneqq\sum_{i=1}^d\frac{1}{\delta_n^2}\int_\Omega R_{\delta_n}(x,y)(x_i-y_i)^2\int_0^1(1-t)\tilde{u}_{ii}(y+t(x-y))dtdy\\
						&-\int_\Omega \bar{R}_{\delta_n}(x,y)\Delta \tilde{u}(y)dy.
					\end{aligned}
				\end{equation*}
				For $r_1^n(x)$, consider the $C^2$ approximation $\tilde{u}_n$ of $\tilde{u}$ satisfying $\left\lVert\tilde{u}-\tilde{u}_n\right\rVert_{H^2(\mathbb{R}^d)}=o(\delta_n^{d+1})$ and $e_n=\tilde{u}-\tilde{u}_n$. According to the Taylor expansion,
				\begin{equation*}
					\begin{aligned}
						0&=\frac{1}{\delta_n^2}\int_\Omega R_{\delta_n}(x,y)(\tilde{u}_n(x)-\tilde{u}_n(y)+\langle\nabla \tilde{u}_n(y),x-y\rangle\\
						&-\sum_{i,j=1}^d\frac{1}{\delta_n^2}\int_\Omega R_{\delta_n}(x,y)(x_i-y_i)(x_j-y_j)\int_0^1(1-t)\tilde{u}_{n,ij}(y+t(x-y))dt)dy
					\end{aligned}
				\end{equation*}
				for all $n$. And the $L^2(\Omega)$ norm of
				\begin{equation*}
					\begin{aligned}
						r_1^n(x)&\coloneqq\frac{1}{\delta_n^2}\int_\Omega R_{\delta_n}(x,y)(e_n(x)-e_n(y)+\langle\nabla e_n(y),x-y\rangle\\
						&-\sum_{i,j=1}^d\frac{1}{\delta_n^2}\int_\Omega R_{\delta_n}(x,y)(x_i-y_i)(x_j-y_j)\int_0^1(1-t)e_{n,ij}(y+t(x-y))dtdy
					\end{aligned}
				\end{equation*}
				tends to $0$ since the $H^2(\mathbb{R}^d)$ norm of $e_n$ decays at a sufficiently rapid rate.\\
				For $r_2^n(x)$, do the decomposition
				\begin{equation*}
					\begin{aligned}
						&\quad\frac{1}{\delta_n^2}\int_\Omega R_{\delta_n}(x,y)(x_i-y_i)(x_j-y_j)\int_0^1(1-t)\tilde{u}_{ij}(y+t(x-y))dtdy\\
						&=\frac{1}{\delta_n^2}\int_\Omega R_{\delta_n}(x,y)(x_i-y_i)(x_j-y_j)\int_0^1(1-t)\tilde{u}_{ij}(x)dtdy\\
							&+\frac{1}{\delta_n^2}\int_\Omega R_{\delta_n}(x,y)(x_i-y_i)(x_j-y_j)\int_0^1(1-t)(\tilde{u}_{ij}(y+t(x-y))-\tilde{u}_{ij}(x))dtdy\\
					\end{aligned}
				\end{equation*}
				For the first term,
				\begin{equation*}
					\begin{aligned}
						&\quad\left\lVert\frac{1}{\delta_n^2}\int_\Omega R_{\delta_n}(x,y)(x_i-y_i)(x_j-y_j)\int_0^1(1-t)\tilde{u}_{ij}(x)dtdy\right\rVert_{L^2(\Omega)}\\
						&=\frac{1}{2}\left\lVert\frac{1}{\delta_n^2}\tilde{u}_{ij}(x)\int_\Omega R_{\delta_n}(x,y)(x_i-y_i)(x_j-y_j)dy\right\rVert_{L^2(\Omega)}\\
						&=\frac{1}{2}\left\lVert\frac{1}{\delta_n^2}\tilde{u}_{ij}(x)\int_\Omega R_{\delta_n}(x,y)(x_i-y_i)(x_j-y_j)dy\right\rVert_{L^2(\{x\in\Omega\big|\mathrm{dist}(x,\partial\Omega)\leq\delta_nr_R\})}\\
						&=O\left(\left\lVert\tilde{u}_{ij}(x)\right\rVert_{L^2\left(\left\{x\in\Omega\big|\mathrm{dist}(x,\partial\Omega)\leq\delta_nr_R\right\}\right)}\right)\\
						&=o(1). \quad\text{as } \delta_n\to 0.
					\end{aligned}
				\end{equation*}
				For the second term,
				\begin{equation*}
					\begin{aligned}
						&\left\lVert\frac{1}{\delta_n^2}\int_\Omega R_{\delta_n}(x,y)(x_i-y_i)(x_j-y_j)\int_0^1(1-t)(\tilde{u}_{ij}(y+t(x-y))-\tilde{u}_{ij}(x))dtdy\right\rVert_{L^2(\Omega)}\\
						&\leq \int_\Omega\left|\int_\Omega R_{\delta_n}(x,y)\int_0^1(\tilde{u}_{ij}(y+t(x-y))-\tilde{u}_{ij}(x))dtdy\right|^2dx\\
						&=O\left(\int_\Omega\int_\Omega R_{\delta_n}(x,y)\int_0^1\left|\tilde{u}_{ij}(y+t(x-y))-\tilde{u}_{ij}(x)\right|^2dtdydx\right)\\
						&=O\left(\int_0^1\int_{\mathbb{R}^d}R(\frac{|z|^2}{4})\int_{\mathbb{R}^d} |\tilde{u}_{ij}(x+(1-t)\delta_n z)-\tilde{u}_{ij}(x)|^2dxdzdt\right)\\
						&=o(1).
					\end{aligned}
				\end{equation*}
				The last estimate is derived from the dominated convergence theorem. Similarly, we have
				\begin{equation*}
					\left\lVert\frac{1}{\delta_n^2}\int_\Omega R_{\delta_n}(x,y)(x_i-y_i)^2\int_0^1(1-t)(\tilde{u}_{i}(y+t(x-y))-\tilde{u}_{ii}(x))dtdy\right\rVert_{L^2(\Omega)}
				\end{equation*}
				tending to $0$ as $n\rightarrow\infty$. Therefore, for $r_3^n(x)$,
				\begin{equation*}
					\begin{aligned}
			&\quad\lim_{n\rightarrow\infty}\left\lVert r_3^n\right\rVert_{L^2(\Omega)}\\
   &\leq\lim_{n\rightarrow\infty}\left\lVert\sum_{i=1}^d\tilde{u}_{ii}(x)\left(\frac{1}{2\delta_n^2}\int_\Omega R_{\delta_n}(x,y)(x_i-y_i)^2dy-\int_\Omega \bar{R}_{\delta_n}(x,y)dy\right)\right\rVert_{L^2(\Omega)}\\
						&=0.
					\end{aligned}
				\end{equation*}
				The last equation is attributed to the fact that
    \begin{equation*}
        \frac{1}{2\delta_n^2}\int_{\mathbb{R}^d} R_{\delta_n}(x,y)(x_i-y_i)^2dy = \frac{\sigma_R}{2} =\int_{\mathbb{R}^d} \bar{R}_{\delta_n}(x,y)dy.
    \end{equation*}
    
    \end{proof}

		Now, we can demonstrate the $\Gamma-$convergence from $F_n$ to $F$. According to the definition, we divide the proof into two parts, corresponding to \cref{lem:liminf} and \cref{lem:limsup} respectively.
		\begin{lemma}[liminf inequality]
			\label{lem:liminf}
			Suppose that $\Omega$ is a bounded domain in $\mathbb{R}^d$ with $C^2$ boundary or a convex polytope. $a\in H^{\frac{3}{2}}(\partial\Omega)$. $b\in H^\frac{1}{2}(\partial \Omega)$.  $K,R,\tilde{R}$ are kernels satisfying {{Assumption} 2.1}. $\{\delta_n\}$, $\{\xi_n\}$ are sequences of positive constants tending to 0. Then, for all $\{u_n\}$ converging to $u$ in $L^2(\Omega)$, we have
			\begin{equation*}
				\liminf_{n\rightarrow\infty}F_n(u_n)\geq F(u),
			\end{equation*}
			where $F_n,F$ are defined as \eqref{eq:nonlocal-functional},\eqref{eq:local-functional}.
		\end{lemma}
		\begin{proof}
			The liminf inequality from $\int_\Omega f u_n dx$ to $\int_\Omega f u dx$ is trivial. Without loss of generality, we can assume that $\liminf_{n\rightarrow\infty}F_n(u_n)<\infty$. Thanks to Theorem 2.1 in \cite{gan2024convergence} and \cref{lem:basic-control}, we have $u\in H^1(\Omega)$, $u=a$  on $\partial\Omega$ and $\liminf_{n\rightarrow\infty}\varepsilon L_n(u_n)+B_n(u_n,a)\geq\varepsilon L(u)$. Subsequently, Note that $b\in H^\frac{1}{2}(\Omega)$. Consider a function $v\in H^2(\Omega)$ which is the solution of the following Laplace equation (existence of $v$ is guaranteed by proposition 7.7 in section 5 of \cite{taylor1996partial}),
			\begin{equation*}
				\left\{\begin{array}{cc}
						-\Delta u\equiv\beta&\text{in }\Omega,\\
						\frac{\partial u}{\partial n}=b&\text{on }\partial\Omega
					\end{array}\right.
			\end{equation*}
			where $\beta=\frac{\int_{\partial\Omega}b(x)dx}{m(\Omega)}$ is a constant. According to \cref{lem:appro}, 
			\begin{equation*}
                    \begin{aligned}
                    \lim_{n\rightarrow\infty}\bigg\lVert\frac{1}{\delta_n^2}\int_{\Omega}R_{\delta_n}(x,y)(v(x)-v(y))dy&-2\int_{\partial\Omega}\bar{R}_{\delta_n}(x,y)b(y)dS(y)\\
                        &-\beta\int_\Omega \bar{R}_{\delta_n}(x,y)dy\bigg\rVert_{L^2(\Omega)}=0.
                    \end{aligned}
			\end{equation*}	
			Hence, for the interior term in $F_n$,
			\begin{equation*}
				\begin{aligned}
					&\quad\liminf_{n\rightarrow\infty}\int_\Omega|\frac{1}{\delta_n^2}\int_{\Omega}R_{\delta_n}(x,y)(u_n(x)-u_n(y))dy-2\int_{\partial\Omega}\bar{R}_{\delta_n}(x,y)b(y)dS(y)|^2dx\\
					&=\liminf_{n\rightarrow\infty}\int_\Omega|\frac{1}{\delta_n^2}\int_{\Omega}R_{\delta_n}(x,y)(w_n(x)-w_n(y))dy+\beta\int_\Omega \bar{R}_{\delta_n}(x,y)dy|^2dx\\
                    &=\liminf_{n\rightarrow\infty}\int_\Omega|\frac{1}{\delta_n^2}\int_{\Omega}R_{\delta_n}(x,y)(w_n(x)-w_n(y))dy+\frac{\beta\sigma_R}{2}|^2dx
				\end{aligned}
			\end{equation*}
			is finite, where $w_n\coloneqq u_n-v$. Since $L^2(\Omega)$ is reflexive, there exists a weak convergence limit $g$ of $\frac{1}{\delta_n^2}\int_{\Omega}R_{\delta_n}(x,y)(w_n(x)-w_n(y))dy+\frac{\beta\sigma_R}{2}$ in the sense of subsequence (still use the original notation for simplicity). This means that for all $\phi\in C^{\infty}(\mathbb{R}^d)$,
			\begin{equation}
				\label{eq:weak-convergence}
				\int_\Omega\phi(x)\frac{1}{\delta_n^2}\int_{\Omega}R_{\delta_n}(x,y)(w_n(x)-w_n(y))dydx\longrightarrow\int_{\Omega}\phi (g-\frac{\beta\sigma_R}{2})dx
			\end{equation}
			as $n\rightarrow\infty$. For the left-hand side,
			\begin{equation*}
				\begin{aligned}
					&\quad\int_\Omega\phi(x)\frac{1}{\delta_n^2}\int_{\Omega}R_{\delta_n}(x,y)(w_n(x)-w_n(y))dydx\\
					&=\int_\Omega\phi(x)\int_{x+\delta_nz\in\Omega}R(\frac{|z|^2}{4})\frac{w_n(x)-w_n(x+\delta_nz)}{\delta_n^2}dzdx\\
					&=\int_\Omega w_n(x)\int_{x\pm\delta_nz\in\Omega} R(\frac{|z|^2}{4})\frac{\phi(x)-\phi(x-\delta_nz)}{\delta_n^2}dzdx\\
					&+\int_\Omega w_n(x)\int_{x+\delta_nz\in\Omega,x-\delta_nz\notin\Omega} R(\frac{|z|^2}{4})\frac{\phi(x)}{\delta_n^2}dzdx\\
					&-\int_\Omega w_n(x)\int_{x-\delta_nz\in\Omega,x+\delta_nz\notin\Omega} R(\frac{|z|^2}{4})\frac{\phi(x-\delta_nz)}{\delta_n^2}dzdx\\
					&=\int_\Omega w_n(x)\int_{x\pm\delta_nz\in\Omega} R(\frac{|z|^2}{4})\frac{\phi(x)-\phi(x-\delta_nz)}{\delta_n^2}dzdx\\
					&+\int_\Omega w_n(x)\int_{x+\delta_nz\in\Omega,x-\delta_nz\notin\Omega} R(\frac{|z|^2}{4})\frac{\phi(x)-\phi(x+\delta_nz)}{\delta_n^2}dzdx\\
					&=-\frac{1}{\delta_n}\int_\Omega w_n(x)\int_{x+\delta_nz\in\Omega} R(\frac{|z|^2}{4})\langle\nabla\phi(x),z\rangle dxdz\\
					&-\frac{1}{2}\int_\Omega w_n(x)\int_{x+\delta_nz\in\Omega} R(\frac{|z|^2}{4})z^TH_\phi(x)zdxdz
					+o(1).
				\end{aligned}
			\end{equation*}
			If the boundary of $\Omega$ is $C^2$, let $\Omega_n\coloneqq\{x\in\Omega\big|\mathrm{dist}(x,\partial\Omega)<\delta_nr_R\}$. Then for large enough $n$, $\Omega_n$ can be reformulated as $\Omega_n=\{x+sn(x)\big|x\in\partial\Omega,s\in(0,\delta_nr_R)\}$ where $n(x)$ is the inward normal vector of $\partial\Omega$ at $x$. Let $w\coloneqq u-v$ For the first term,
			\begin{equation*}
				\begin{aligned}
					&\quad\lim_{n\rightarrow\infty}\frac{-1}{\delta_n}\int_\Omega w_n(x)\int_{x+\delta_nz\in\Omega} R(\frac{|z|^2}{4})\langle\nabla\phi(x),z\rangle dzdx\\
					&=\lim_{n\rightarrow\infty}\frac{-1}{\delta_n}\int_{\Omega_n} w(x)\int_{x+\delta_nz\in\Omega}R(\frac{|z|^2}{4})\langle\nabla\phi(x),z\rangle dzdx\\
					&=\lim_{n\rightarrow\infty}\frac{-1}{\delta_n}\int_{\partial\Omega}w(x)\frac{\partial\phi}{\partial n}(x)\int_0^{\delta_nr_R}\int_{x+sn(x)+\delta_nz\in\Omega}R(\frac{|z|^2}{4})\langle n(x),z\rangle dzdsdS(x)\\
					&=\lim_{n\rightarrow\infty}\frac{-1}{\delta_n}\int_{\partial\Omega}w(x)\frac{\partial\phi}{\partial n}(x)\int_{\langle n(x),z\rangle\geq 0}\int_0^{\delta_nr_R}R(\frac{|z|^2}{4})\langle n(x),z\rangle dsdzdS(x)\\
					&-\lim_{n\rightarrow\infty}\frac{1}{\delta_n}\int_{\partial\Omega}w(x)\frac{\partial\phi}{\partial n}(x)\int_{\langle n(x),z\rangle<0}\int_{-\delta_n \langle n(x),z\rangle}^{\delta_nr_R}R(\frac{|z|^2}{4})\langle n(x),z\rangle dsdzdS(x)\\
					&=\int_{\partial\Omega}w(x)\frac{\partial\phi}{\partial n}(x)\int_{\langle n(x),z\rangle\geq 0}R(\frac{|z|^2}{4})\langle n(x),z\rangle^2dzdx\\
					&=\frac{1}{2}\int_{\mathbb{R}^d}R(\frac{|z|^2}{4})z_1^2dz\int_{\partial\Omega}w(x)\frac{\partial\phi}{\partial n}(x)dx.
				\end{aligned}				
			\end{equation*}
            If $\Omega$ is a polytope, we denote $$\partial\Omega = \bigcup_{i=1}^N\mathcal{F}_i\cup\mathcal{R}$$ where $\left\{\mathcal{F}_i\right\}_{i=1}^N$ are the facets and $\mathcal{R}$ is the ridge. Then, there exists a constant $C>0$ only depend on $\Omega$ and continuous functions $\alpha_n(x):\partial\Omega\rightarrow[0,1]$ satisfying $\alpha_n(x)=1$ for all $x$ fulfilling $\mbox{dist}(x,\mathcal{R})>C\delta_n$ for large enough $n$ such that $\Omega_n$ can be reformulated as $$\Omega_n=\{x+s\alpha_n(x)n(x)\big|x\in\partial\Omega,s\in(0,\delta_nr_R)\}.$$ In fact, we can define 
            $$\Omega_n^i\coloneqq\{x\in\Omega_n\big|\mbox{dist}(x,\mathcal{F}_i)=\min_{1\leq j\leq N}\mbox{dist}(x,\mathcal{F}_j)\}$$ and $$\alpha_n(x)\coloneqq\sup\{\alpha\in[0,1]\big|x+\delta_nr_R\alpha n(x)\in\Omega_n^i\}\quad \mbox{for all}\quad x\in\mathcal{F}_i.$$ Note that $\alpha_n(x)\rightarrow 1$ as $n\rightarrow\infty$ a.s. The equation above still holds because of the Lebesgue dominated convergence theorem.
            \begin{equation*}
				\begin{aligned}					&\quad\lim_{n\rightarrow\infty}\frac{-1}{\delta_n}\int_\Omega w_n(x)\int_{x+\delta_nz\in\Omega} R(\frac{|z|^2}{4})\langle\nabla\phi(x),z\rangle dzdx\\
					&=\lim_{n\rightarrow\infty}\frac{-1}{\delta_n}\int_{\partial\Omega}w(x)\frac{\partial\phi}{\partial n}(x)\int_0^{\delta_nr_R\alpha_n(x)}\int_{x+sn(x)+\delta_nz\in\Omega}R(\frac{|z|^2}{4})\langle n(x),z\rangle dzdsdS(x)\\
					&=\lim_{n\rightarrow\infty}\frac{-1}{\delta_n}\int_{\partial\Omega}w(x)\frac{\partial\phi}{\partial n}(x)\int_0^{\delta_nr_R}\int_{x+sn(x)+\delta_nz\in\Omega}R(\frac{|z|^2}{4})\langle n(x),z\rangle dzdsdS(x)\\
					&=\frac{1}{2}\int_{\mathbb{R}^d}R(\frac{|z|^2}{4})z_1^2dz\int_{\partial\Omega}w(x)\frac{\partial\phi}{\partial n}(x)dx.
				\end{aligned}				
			\end{equation*}
			For the second term,
			\begin{equation*}
				\begin{aligned}
					&\quad\lim_{n\rightarrow\infty}-\frac{1}{2}\int_\Omega w_n(x)\int_{x+\delta_nz\in\Omega} R(\frac{|z|^2}{4})z^TH_\phi(x)zdxdz\\
					&=\lim_{n\rightarrow\infty}-\frac{1}{2}\int_\Omega w(x)\int_{x+\delta_nz\in\Omega} R(\frac{|z|^2}{4})z^TH_\phi(x)zdxdz\\
					&=\lim_{n\rightarrow\infty}-\frac{1}{2}\int_{\Omega\backslash\Omega_n} w(x)\int_{x+\delta_nz\in\Omega} R(\frac{|z|^2}{4})z^TH_\phi(x)zdxdz\\
					&=\lim_{n\rightarrow\infty}-\frac{1}{2}\int_{\Omega\backslash\Omega_n} w(x)\int_{\mathbb{R}^d} R(\frac{|z|^2}{4})z^TH_\phi(x)zdxdz\\
					&=-\frac{1}{2}\sum_{i=1}^d \int_\Omega w(x)\int_{\mathbb{R}^d} R(\frac{|z|^2}{4})\phi_{ii}(x)z_i^2dxdz\\
					&=-\frac{1}{2}\int_{\mathbb{R}^d}R(\frac{|z|^2}{4})z_1^2dz\int_\Omega w(x)\Delta\phi(x)dx.
				\end{aligned}
			\end{equation*}
			Therefore, we have
			\begin{equation*}
\left(\int_{\partial\Omega}w(x)\frac{\partial\phi}{\partial n}(x)dx-\int_\Omega w(x)\Delta\phi(x)dx\right)=\int_\Omega\phi (g-\frac{\beta\sigma_R}{2}) dx.
			\end{equation*} 
			It means that
			\begin{equation*}
				\int_\Omega \nabla w(x)\cdot\nabla\phi(x)dx=\int\phi ( 
                g-\beta)dx
			\end{equation*}
			and $w$ is the weak solution of 
			\begin{equation*}
				\left\{\begin{array}{cc}
					\Delta w= 
                g-\beta&\text{in }\Omega,\\
					w=a-v&\text{on }\partial\Omega.
				\end{array}\right.
			\end{equation*}
			Note that $a\in H^\frac{3}{2}(\partial\Omega)$. According to the regularity of weak solutions of elliptic equation (see, for example, in proposition 1.7 in section 5 of \cite{taylor1996partial} and Theorem 3.2.1.2 of \cite{grisvard2011elliptic}), we have $w\in H^2(\Omega)$. From \eqref{eq:weak-convergence}, it also yields that
			\begin{equation*}
				\begin{aligned}
					&\quad\liminf_{n\rightarrow\infty}\int_\Omega\left|\frac{1}{\delta_n^2}\int_{\Omega}R_{\delta_n}(x,y)(u_n(x)-u_n(y))dy-2\int_{\partial\Omega}\bar{R}_{\delta_n}(x,y)b(y)dS(y)\right|^2dx\\
					&=\liminf_{n\rightarrow\infty}\int_\Omega\left|\frac{1}{\delta_n^2}\int_{\Omega}R_{\delta_n}(x,y)(w_n(x)-w_n(y))dy+\frac{\beta\sigma_R}{2}\right|^2dx\\
					&\geq \int_\Omega g^2 dx=\frac{\sigma_R^2}{4}\int_{\Omega}|\Delta w+\beta|^2dx=\frac{\sigma_R^2}{4}\int_{\Omega}|\Delta u|^2dx.
				\end{aligned}
			\end{equation*}
			Finally, we only need to prove that $\frac{\partial u}{\partial n}=b$ with
			\begin{equation*}
				\liminf_{n\rightarrow\infty}\int_\Omega\left|\frac{1}{\delta_n^2}\int_{\Omega}R_{\delta_n}(x,y)(u_n(x)-u_n(y))dy-2\int_{\partial\Omega}\bar{R}_{\delta_n}(x,y)b(y)dS(y)\right|^2dx<\infty.
			\end{equation*}
			It means that
			\begin{equation*}
				\liminf_{n\rightarrow\infty}\int_\Omega\left|\frac{1}{\delta_n^2}\int_{\Omega}R_{\delta_n}(x,y)(u(x)-u(y))dy-2\int_{\partial\Omega}\bar{R}_{\delta_n}(x,y)b(y)dS(y)\right|^2dx<\infty.
			\end{equation*}
			And thanks to \cref{lem:appro}, we have
			\begin{equation*}
				\lim_{n\rightarrow\infty}\int_\Omega\left|\int_{\partial\Omega}\bar{R}_{\delta_n}(x,y)(\frac{\partial u}{\partial n}(y)-b(y))dS(y)\right|^2dx<\infty
			\end{equation*}
			and
			\begin{equation*}
				\lim_{n\rightarrow\infty}\delta_n\int_\Omega\left|\int_{\partial\Omega}\bar{R}_{\delta_n}(x,y)(\frac{\partial u}{\partial n}(y)-b(y))dS(y)\right|^2dx=0
			\end{equation*}
			For large enough $n$,
			\begin{equation*}
				\begin{aligned}
					&\quad\delta_n\int_\Omega\left|\int_{\partial\Omega}\bar{R}_{\delta_n}(x,y)(\frac{\partial u}{\partial n}(y)-b(y))dS(y)\right|^2dx\\
					&=\delta_n\int_{\Omega_n}\left|\int_{\partial\Omega}\bar{R}_{\delta_n}(x,y)(\frac{\partial u}{\partial n}(y)-b(y))dS(y)\right|^2dx\\
					&=\delta_n\int_{\partial\Omega}\int_0^{\delta_n r_{\bar{R}}}\left|\int_{\partial\Omega}\bar{R}_{\delta_n}(x+sn(x),y)(\frac{\partial u}{\partial n}(y)-b(y))dS(y)\right|^2dsdS(x)\\
					&\geq \frac{1}{r_{\bar{R}}}\int_{\partial\Omega}\left|\int_0^{\delta_n r_{\bar{R}}}\int_{\partial\Omega}\bar{R}_{\delta_n}(x+sn(x),y)(\frac{\partial u}{\partial n}(y)-b(y))dS(y)ds\right|^2dS(x)\\
					&=\frac{1}{r_{\bar{R}}}\int_{\partial\Omega}\left|\int_{\partial\Omega}\bar{\bar{R}}_{\delta_n}(x,y)(\frac{\partial u}{\partial n}(y)-b(y))dS(y)\right|^2dS(x)
				\end{aligned}
			\end{equation*}
			where $\bar{\bar{R}}_{\delta_n}(x,y)\coloneqq\int_0^{\delta_n r_{\bar{R}}}\bar{R}_{\delta_n}(x+sn(x),y)ds=\frac{1}{\delta_n^{d-1}}\int_0^{ r_{\bar{R}}}\bar{R}(\frac{|\frac{x-y}{\delta_n}+sn(x)|^2}{4})ds.$ Hence,
			\begin{equation*}
				\begin{aligned}
					&\quad\lim_{n\rightarrow\infty}\frac{1}{r_{\bar{R}}}\int_{\partial\Omega}\left|\int_{\partial\Omega}\bar{\bar{R}}_{\delta_n}(x,y)(\frac{\partial u}{\partial n}(y)-b(y))dS(y)\right|^2dS(x)\\
					&=\frac{1}{r_{\bar{R}}}\int_{\partial\Omega}|\frac{\partial u}{\partial n}(x)-b(x)|^2dS(x)\\
				\end{aligned}
			\end{equation*}
			which means that $\left\lVert\frac{\partial u}{\partial n}(x)-b(x)\right\rVert_{L^2(\partial\Omega)}=0$ and $\frac{\partial u}{\partial n}=b$ almost surely on $\partial\Omega$.
		\end{proof}
		
		\begin{lemma}[limsup inequality]
			\label{lem:limsup}
			Suppose that $\Omega$ is a bounded domain in $\mathbb{R}^d$ with $C^2$ boundary.  $K,R,\tilde{R}$ are kernels satisfying {{Assumption} 2.1}. $\{\delta_n\}$, $\{\xi_n\}$ are sequences of positive constants tending to 0 satisfying $\lim_{n\rightarrow\infty}\frac{\delta_n}{\xi_n}=0$. Then, for all $u$ in $L^2(\Omega)$, there exists a sequence $\{u_n\}$ converging to $u$ in $L^2(\Omega)$ satisfying
			\begin{equation*}
				\limsup_{n\rightarrow\infty}F_n(u_n)\leq F(u),
			\end{equation*}
			where $F_n,F$ are defined as \eqref{eq:nonlocal-functional},\eqref{eq:local-functional}.
		\end{lemma}
		\begin{proof}
			Without loss of generality, we can assume that $F(u)<\infty$. With such an assumption, we have $u\in H^2(\Omega)$, $u=a$ and $\frac{\partial u}{\partial n}=b$ on $\partial\Omega$. Select $\{u_n\}\in C^\infty(\Omega)$ satisfying $\left\lVert u_n- u\right\rVert_{H^2(\Omega)}=o(\delta_n^\frac{d+2}{2})$. Then, $\lim_{n\rightarrow\infty}\int_\Omega fu_ndx=\int_\Omega fudx$. $\limsup_{n\rightarrow\infty}L_n(u_n)\leq L(u)$ can be derived from Lemma 5.8 in \cite{gan2024convergence}. For the Dirichlet boundary term,
			\begin{equation*}
				\begin{aligned}
					B_n(u_n,a)&\leq\frac{1}{\xi_n}\int_{\partial\Omega}\left|\int_\Omega K_{\delta_n}(|x-y|)(a(x)-u_n(x))dy\right|^2dx\\
					&+\frac{1}{\xi_n}\int_{\partial\Omega}\left|\int_\Omega K_{\delta_n}(|x-y|)(u_n(x)-u_n(y))dy\right|^2dx\\
					&=O\left(\frac{1}{\xi_n}\int_{\partial\Omega}\int_\Omega K_{\delta_n}(|x-y|)|a(x)-u_n(x)|^2dydx\right)\\
					&+O\left(\frac{1}{\xi_n}\int_{\partial\Omega}\int_\Omega K_{\delta_n}(|x-y|)|x-y|^2dydx\right)\\
					&=O\left(\frac{\left\lVert u_n-a\right\rVert_{L^2(\partial\Omega)}^2}{\xi_n}\right)+O\left(\frac{\delta_n^2}{\xi_n}\right)\\
					&=o(1), \quad \text{as } \delta_n\to 0.
				\end{aligned}
			\end{equation*}
			The last control holds owing to the trace theorem of Sobolev functions and $\lim_{n\rightarrow\infty}\frac{\delta_n^2}{\xi_n}\\
            =0$. Finally, for the interior term, applying \cref{lem:appro} again, we have
			\begin{equation*}
				\begin{aligned}
					&\quad\limsup_{n\rightarrow\infty}\int_\Omega\left|\frac{1}{\delta_n^2}\int_{\Omega}R_{\delta_n}(x,y)(u_n(x)-u_n(y))dy-2\int_{\partial\Omega}\bar{R}_{\delta_n}(x,y)b(y)dS(y)\right|^2dx\\
				&=\limsup_{n\rightarrow\infty}\int_\Omega\left|\frac{1}{\delta_n^2}\int_{\Omega}R_{\delta_n}(x,y)(u(x)-u(y))dy-2\int_{\partial\Omega}\bar{R}_{\delta_n}(x,y)b(y)dS(y)\right|^2dx\\
					&=\lim_{n\rightarrow\infty}\int_\Omega\left|\int_\Omega \bar{R}_{\delta_n}(x,y)\Delta u(y)dy\right|^2dx\\
					&= 
                    \int_\Omega|\Delta u|^2dx.
				\end{aligned}
			\end{equation*}
            The first equation is also attributed to the sufficiently rapid rate for the convergence of $u_n$ to $u$.
            \end{proof}
\section{Convergence of minimizers} The convergence of minimizers is rendered in this section. This is actually a key property of $\Gamma-$convergence, which is also applied in \cite{garcia2016continuum,slepcev2019analysis,roith2022continuum,gan2022non}.
We include here for completeness and easy reference.
\begin{lemma}
		\label{lem:minimizer}
		Let $X$ be a metric space and $F_n:X\rightarrow[0,\infty]$ $\varGamma$-converges to $F:X\rightarrow[0,\infty]$ which is not identically $\infty$. If there exists a relatively compact sequence $\{x_n\}_{n\in\mathbb{N}}\subset X$ such that
		\[\lim_{n\rightarrow\infty}(F_n(x_n)-\inf_{x\in X}F_n(x))=0\]
		then we have
		\[\lim_{n\rightarrow\infty}\inf_{x\in X}F_n(x)=\min_{x\in X}F(x)\]
		and any cluster point of $\{x_n\}_{n\in\mathbb{N}}$ is a minimizer of $F$.
\end{lemma}

Having established the $\Gamma-$convergence of $F_n$ in \cref{thm:gamma-convergence}, it remains to demonstrate the compactness to prove the convergence of the minimizers according to \cref{lem:minimizer}. And the following conclusion presented in \cite{gan2024convergence} is helpful.

\begin{lemma}[Lemma 6.1 in \cite{gan2024convergence}]
     \label{lem:compactness}
		Suppose that $\Omega$ is a Lipschitz bounded domain in $\mathbb{R}^d$. $\bar{R}$ is a kernel satisfying {{Assumption} 2.1}. $\{\delta_n\}$ is a sequence of positive constants tending to $0$ as $n\rightarrow\infty$. $\{u_n\}$ is a bounded sequence in $L^2(\Omega)$ and satisfies
		\[\sup_nL_n(u_n)<\infty\]
		where $L_n$ is defined as \eqref{eq:nonlocol-laplacian}. Then $\{u_n\}$ is a relatively compact sequence in $L^2(\Omega)$.
\end{lemma}

\begin{proof}[proof of \cref{thm:convergence-minimizers}]
    With \cref{lem:poincare}, any sequence $\{u_n\}$ satisfying \\$\sup_{n}F_n(u_n)<+\infty$ is bounded in $L^2(\Omega)$. Therefore, \cref{thm:convergence-minimizers} can be directly derived with \cref{lem:minimizer}, \cref{lem:compactness} and \cref{lem:basic-control}.
\end{proof}

In practice, we usually solve the minimizer of nonlocal model $F_n$ numerically. It is necessary to discretize the integrals and approximate $F_n$ with some discrete functionals. The following theorem ensures the convergence of minimizers if the truncation error can be controlled. We required the boundedness of the minimizers of $F_n$. This is guaranteed with \cref{lem:boundedness-minimizers} for the homogeneous boundary condition.

\begin{theorem}

\label{thm:approximate-error}
Let $\{F_n\}_{n=1}^\infty$ be a sequence of functionals on $L^2(\Omega)$ and $B$ be a bounded subset of $L^2(\Omega)$ containing all minimizers of $\{F_n\}_{n=1}^\infty$. Let $\{\tilde{F}_n\}_{n=1}^\infty$ be a sequence of functionals satisfying
\begin{equation*}
    |\tilde{F_n}(u)-F_n(u)|\leq \sigma_n\lVert u\rVert_{L^2(\Omega)}
\end{equation*}
for all $u\in B$ and some constants $\sigma_n\rightarrow 0$. Let $\{v_n\}_{n=1}^\infty\subset B$ be a minimizing sequence of $\{\tilde{F}_n\}_{n=1}^\infty$ which means that 
\begin{equation*}
    \lim_{n\rightarrow\infty}(\tilde{F}_n(v_n)-\inf_{u}\tilde{F}_n(u))=0.
\end{equation*}
Then, 
\begin{equation*}
    \lim_{n\rightarrow\infty}\inf_u\tilde{F}_n(u)=\min_u F(u)
\end{equation*}
and any cluster point of $\{v_n\}_{n=1}^\infty$ with respect to $\lVert\cdot\rVert_{L^2(\Omega)}$ is the minimizer of $F$ in $L^2(\Omega)$.
\end{theorem}

\begin{proof}
    With \cref{thm:convergence-minimizers},
    \begin{equation*}
    \lim_{n\rightarrow\infty}\inf_uF_n(u)=\min_u F(u).
    \end{equation*}
    With \cref{thm:well-posedness}, $\inf_{u}F_n(u)=F_n(u_n^*)>-\infty$ for some unique $u_n^*$. And with \cref{lem:well-posedness2}, there exists a constant $C$ not depending on $n$ such that
    \begin{equation*}
        F_n(v_n)-F_n(u_n^*)\geq C\left\lVert v_{n}-u_{n}^*\right\rVert_{L^2(\Omega)}^2.
    \end{equation*}
    Let $v$ be a cluster point of $\{v_n\}_{n=1}^\infty$. Since $\{u_n^*\}$ is relatively compact, we can select a subsequence satisfying $v_{n_k}-u_{n_k}^*\rightarrow v-u^*$ in $L^2(\Omega)$ where $u^*$ is the minimizer of $F$. Then,
    \begin{equation*}
        \limsup_{k\rightarrow\infty}(F_n(v_{n_k})-F_n(u_{n_k}^*))\geq C\left\lVert v_{n_k}-u_{n_k}^*\right\rVert_{L^2(\Omega)}.
    \end{equation*}
    Define $M=\sup_{u\in B}\lVert u\rVert_{L^2(\Omega)}$. Note that
    \begin{equation*}
        |F_n(v_n)-F_n(u_n^*)|\leq |\tilde{F}_n(v_n)-\inf_{u}\tilde{F}_n(u)|+2\sigma_nM
    \end{equation*}
    and the condition
    \begin{equation*}
    \lim_{n\rightarrow\infty}(\tilde{F}_n(v_n)-\inf_{u}\tilde{F}_n(u))=0.
\end{equation*}
We have 
\begin{equation*}
\lim_{k\rightarrow\infty}\left\lVert v_{n_k}-u_{n_k}^*\right\rVert_{L^2(\Omega)}=0
\end{equation*}
which means that $v=u$ with respect to the $L^2(\Omega)$ norm.
\end{proof}

\section{Numerical Validation}
    In this section, we will verify the convergence of the nonlocal biharmonic model~\cref{eq:nonlocal-functional} to its local counterpart using numerical examples.

    For simplicity of notation, we denote
    \[g(x)=\frac{1}{\delta_n^2}\int_{\Omega}R_{\delta_n}(x,y)(u(x)-u(y))dy-2\int_{\partial\Omega}\bar{R}_{\delta_n}(x,y)b(y)dS(y),\]
    \[h(x)=\int_\Omega K_{\delta_n}(x,y)(a(x)-u(y))d y,\]
    and the nonlocal model~\cref{eq:nonlocal-functional} becomes
    \[F_n(u)=\int_{\Omega} |g(x)|^2 dx -2 \int_\Omega f(x)u(x)dx +\frac{1}{\xi_n} \int_{\partial \Omega} |h(x)|^2 dx.\]
    For a general kernel function $R$ and a computational domain $\Omega\subset \mathbb{R}^d$, discretizing the above functional is not trivial given the computational cost of evaluating the numerical quadratures of dimension $2d$. 

    To simplify the computation, we choose the kernel function as
    $$R_\delta(x,y) =K_\delta(x,y) = c_\delta \exp \left(-\frac{|x-y|^2}{\delta^2} \right) $$ 
    for a normalization constant $c_\delta=4\pi^{-d/2}\delta^{-d}$,
    and the computational domain is taken to be the unit cube $\Omega=[0,1]^d$. The advantage of the above setting is that the multidimensional integrals can be separated in dimension, which reduces the computational cost significantly.    

\subsection{Numerical method}

 We use the finite element method to solve the variational problem \eqref{eq:nonlocal-functional}, i.e. we solve the approximate solution in a finite dimensional space $V_h$, 
$$\min_{u\in V_h} F_n(u),$$
and $V_h$ is constructed by piecewise cubic polynomial basis. 

    First, we use the tensor product mesh
    $$\Omega=\bigcup_{i_1,\cdots,i_d=0}^{N-1} \Omega_{i_1,\cdots,i_d}$$
    with $$\Omega_{i_1,\cdots,i_d}=e_{i_1}\times \cdots \times e_{i_d},\quad e_i=[ih,(i+1)h],\quad h=1/N.$$ 
    
    {The finite element space $V_h$ is given as
    \begin{equation}
    \label{eq:FEM-space}
    V_h=\{u\in C(\Omega): u\; \mbox{is multi-cubic polynomial in} \;\Omega_{i_1,\cdots,i_d}\}.
    \end{equation}
 For the local biharmonic equation, the conformal finite element method requires that the finite element space belongs to $H^2(\Omega)$. To ensure that the nonlocal solution converges to the local solution as $\delta\rightarrow 0$, $V_h\cap H^2(\Omega)$ should be large enough to approximate the local solution. $V_h$ in \eqref{eq:FEM-space} fits this requirement since 
 $$\{u\in C^1(\Omega): u\; \mbox{is multi-cubic polynomial in} \;\Omega_{i_1,\cdots,i_d}\} \subset V_h\cap H^2(\Omega).$$
 Here we use cubic polynomials based on the trade-off between accuracy and complexity of implementation. It is straightforward to construct $V_h$ using higher-order polynomials.}
 
    We use the standard multi-cubic polynomial basis, for $x=(x_1,\cdots,x_d)$,   $$\psi_{j_1,\cdots,j_d}(x)=\psi_{j_1}(x_1)\cdots\psi_{j_d}(x_d),\quad j_1,\cdots,j_d=0,1,\cdots,3N.$$
    Along each dimension, $\psi_j(s)$ is a piecewise cubic polynomial defined as follows
    \[\psi_{3i}(s)=\begin{cases}
        \frac{(s-s_{3i-3})(s-s_{3i-2})(s-s_{3i-1})}{(s_{3i}-s_{3i-3})(s_{3i}-s_{3i-2})(s_{3i}-s_{3i-1})}, & \quad s\in e_{i-1} \\
        \frac{(s-s_{3i+1})(s-s_{3i+2})(s-s_{3i+3})}{(s_{3i}-s_{3i+1})(s_{3i}-s_{3i+2})(s_{3i}-s_{3i+3})}, & \quad s\in e_{i} \\
        0, & \quad \text{else}
    \end{cases}\]
    \[\psi_{3i+1}(s)=\begin{cases}
        \frac{(s-s_{3i})(s-s_{3i+2})(s-s_{3i+3})}{(s_{3i+1}-s_{3i})(s_{3i+1}-s_{3i+2})(s_{3i+1}-s_{3i+3})}, & \quad s\in e_{i} \\
        0, & \quad \text{else}
    \end{cases}\]
    \[\psi_{3i+2}(s)=\begin{cases}
        \frac{(s-s_{3i})(s-s_{3i+1})(s-s_{3i+3})}{(s_{3i+2}-s_{3i})(s_{3i+2}-s_{3i+1})(s_{3i+2}-s_{3i+3})}, & \quad s\in e_{i} \\
        0, & \quad \text{else}
    \end{cases}\]
    and $s_j=jh/3$.
   { Along with the Gaussian kernel, for the basis with tensor product form, it is easy to transform the $2d$ dimensional integral to product of 2 dimensional integrals, which reduces the computational cost and complexity tremendously.}
    
    For any $u\in V_h$, we need to compute $F_n(u)$. Let    $$u(x)=\sum_{j_1,\cdots,j_d=0}^{3N} u_{j_1,\cdots,j_d} \psi_{j_1,\cdots,j_d}(x).$$

    First, we focus on the inner integral $g(x)$ and $h(x)$ in the nonlocal model. To simplify the notation, we take $d=2$ and take the integral $\int_\Omega R_{\delta_n}(x,y) u(y) dy$ as an example.
    \begin{align*}
        & \int_\Omega R_{\delta_n}(x,y) u(y) dy  =\sum_{j_1,j_2=0}^{3N}\int_{\Omega} R_{\delta_n}(x,y) u_{j_1,j_2}\psi_{j_1,j_2}(y) dy \\
        =& \sum_{i_1,i_2=0}^{N-1}\sum_{j_1,j_2=0}^{3N}\int_{i_1h}^{(i_1+1)h}\int_{i_2h}^{(i_2+1)h} R_{\delta_n}(x,y) u_{j_1,j_2}\psi_{j_1,j_2}(y) dy \\
        ={} & \sum_{i_1,i_2=0}^{N-1}\sum_{j_1,j_2=0}^{3N} c_\delta \cdot \text{coeff}(x,i_1,i_2,j_1,j_2) \cdot u_{j_1,j_2}
    \end{align*}
   The coefficient $\text{coeff}(x,i_1,i_2,j_1,j_2)$ is defined as
    \begin{align*}
        & \text{coeff}(x,i_1,i_2,j_1,j_2) \\
        ={} & \int_{i_1h}^{i_1h+h} 
        \int_{i_2h}^{i_2h+h} 
        \exp\left(-\frac{(x_1-y_1)^2+(x_2-y_2)^2}{\delta_n^2}\right) \psi_{j_1}(y_1)\psi_{j_2}(y_2) \, dy_2 \, dy_1 \\
        ={} & \int_{i_1h}^{i_1h+h} \exp\left(-\frac{(x_1-y_1)^2}{\delta_n^2}\right) \psi_{j_1}(y_1)  d y_1 \cdot \int_{i_2h}^{i_2h+h} \exp\left(-\frac{(x_2-y_2)^2}{\delta_n^2}\right) \psi_{j_2}(y_2)  d y_2 
    \end{align*}
    In each cell $[ih,ih+h]$, the basis function is cubic polynomial and the integrals above can be computed explicitly. Let
    \[ f_0(\eta,a,b)=\int_a^b \exp (-\eta^2r^2) dr= \frac{\sqrt{\pi}}{2\eta} \left(\text{erf}(\eta b) - \text{erf}(\eta a) \right) \]
	\[ f_1(\eta ,a,b)=\int_a^b r \exp (-\eta ^2r^2) dr= \frac{1}{2\eta ^2} (\exp(-\eta ^2a^2)-\exp(-\eta ^2b^2)) \]
	\[ f_2(\eta ,a,b)=\int_a^b r^2 \exp (-\eta ^2r^2) dr= \frac{a}{2\eta ^2}\exp(-\eta ^2a^2) - \frac{b}{2\eta ^2}\exp(-\eta ^2b^2) + \frac{1}{2\eta ^2}f_0(\eta ,a,b)\]
	\[ f_3(\eta ,a,b)=\int_a^b r^3 \exp (-\eta ^2r^2) dr= \frac{a^2}{2\eta ^2}\exp(-\eta ^2a^2) - \frac{b^2}{2\eta ^2}\exp(-\eta ^2b^2) + \frac{1}{\eta ^2} f_1(\eta ,a,b)\]
    where $\text{erf}(z)=\frac{2}{\sqrt{\pi}}\int_0^z e^{-t^2} \mathrm{d} t$ represents Gauss error function. $\eta$ is the scaling factor, which equals $1/\delta_n$ in our experiments. Given that each basis function $\psi_j$ is a cubic polynomial, $\text{coeff}(x,i_1,i_2,j_1,j_2)$ can be represented by linear combination of $f_0,f_1,f_2$ and $f_3$ and the coefficients of the linear combination are determined by the coefficients of the cubic polynomials. 


    The discretization of the boundary integral terms is treated similarly. If $d=1$, then $\partial \Omega=\{0,1\}$ is a two-point boundary, and the boundary integral becomes the sum of the function values at these two points. If $d=2$, then $\partial \Omega$ is the union of four one-dimensional line segments, and we can use the same quadrature rule presented earlier to perform the discretization. Thus far, we are able to evaluate the function value of $g(x)$ and $h(x)$ at any point using $u_{j_1,\cdots,j_d}$, the function value on basis points.

    Next, we present the discretization method for the outer integral. We focus on the one-dimensional case and provide the integration points $x_i$ and quadrature weights $w_i$. The extension to the two-dimensional case is straightforward: the integration points are chosen as the Cartesian product of the one-dimensional points in each direction, and the corresponding quadrature weight at each point is given by $w_{ij}^{\text{2D}} =w_{i}^{\text{1D}} \cdot w_{j}^{\text{1D}}$, which can be efficiently implemented using the Kronecker product.
 
    The term $\int_\Omega f(x)u(x)dx$ is computed using the composite Simpson's 3/8 rule on each sub-interval $e_i=[ih,(i+1)h]$ , as this term exhibits sufficient regularity. In contrast, the term $\int_{\Omega} |g(x)|^2 dx$ requires a different treatment, since our choice of piecewise cubic polynomial basis functions does not belong to $C^2$. As a result, the term $\frac{1}{\delta_n^2}\int_{\Omega}R_{\delta_n}(x,y)(u(x)-u(y))dy$ in $g(x)$, which approximates $\Delta u(x)$, exhibits fast change in a layer near the boundaries of each sub-interval $e_i$, and the width of the layer is $O(\delta)$. This may lead to large error of the outer integral if equally spaced integration points are used, as in Simpson’s rule. 
    
    To address this issue, we propose to handle the boundary layers and interior area separately as following:
    \begin{itemize}[leftmargin=*]
        \item If $2\times 3\delta < h$, we apply the Gauss-Legendre quadrature rule with 5 points on each of the three sub-intervals: $[ih,ih+3\delta], [ih+3\delta, (i+1)h-3\delta], [(i+1)h-3\delta,(i+1)h]$.
        \item If $2\times 3\delta \geq h$, we apply the Gauss-Legendre quadrature rule with 15 points over the entire sub-interval $e_i=[ih,(i+1)h]$.
    \end{itemize}
    In this method, we treat the regions of length $3\delta$ near the interval boundaries as a  boundary layer, if $\delta$ is small. The choice of $3\delta$ is motivated by the empirical $3\sigma$-rule for the Gaussian distribution~\cite{Kazmier2009schaum}. If $\delta$ is sufficiently large, then the entire sub-interval is regarded as an interior region. The Gauss-Legendre quadrature rule is originally defined on the reference interval $[-1,1]$, with integration points $\tilde{x}_i$ and weights $\tilde{w}_i$. To apply the rule to any interval $[a,b]$, we rescale the points and weights as
    \[x_i=\frac{b-a}{2} \tilde{x}_i+\frac{a+b}{2}, \quad w_i=\frac{b-a}{2} \tilde{w}_i\]
    The final set of integration points on $\Omega=[0,1]$ is obtained by collecting the integration points from each sub-interval $e_i$.
    
    The last term, $\int_{\partial \Omega} |h(x)|^2 dx$, is trivial in the one-dimensional case, since $\partial \Omega=\{0,1\}$. In the two-dimensional case, we use the same method as above to compute the integral of $|g(x)|^2$ over each of the four boundary line segments.

    Finally, we will verify that the minimizer of the nonlocal functional $F_n(u)$ converges to the minimizer of its local counterpart $F(u)$, which is the solution of the biharmonic equation $\Delta^2u=f$. With the discretization methods described above, $F_n(u)$ can be written as a quadratic form in terms of the vector $\mathbf{u}=[u_{j_1}]_{j_1=0}^{3N} \in \mathbb{R}^{3N+1}$ for $d=1$, and $\mathbf{u}=[u_{j_1,j_2}]_{j_1,j_2=0}^{3N} \in \mathbb{R}^{(3N+1)^2}$ for $d=2$. Then, taking the derivative with respect to $\mathbf{u}$, we obtain a linear system. Thanks to the quadratic structure of the functional, the resulting linear system is symmetric and positive definite. The numerical experiments in this paper is to verify the theoretical results, so we do not conduct large scale experiments. The size of the linear system is moderate and solved by direct method. 
    For systems of a much larger size, iterative methods may be more efficient. We will carry out more comprehensive studies on the numerical solver in the future work.
    
    \subsection{Numerical results}
    We conduct numerical experiments in 1D and 2D to check the convergence of the solution of the nonlocal model to the solution of the local biharmonic model. For 1D example, we choose the exact solution of the local biharmonic model as $u_{\text{gt}}(x)=x^{10}$. For 2D example, the exact solution of the local model is $u_{\text{gt}}(x_1,x_2)=x_1 \ln (x_2+1)$. Then $f$ is analytically computed by the local biharmonic equation $f(\mathbf{x})=\Delta^2 u(\mathbf{x})$. The boundary data are included in the nonlocal model in the form of $a(\mathbf{x})=u_{\text{gt}}(\mathbf{x})$ and $b(\mathbf{x})=\frac{\partial u_{\text{gt}}}{\partial \mathbf{n}}(\mathbf{x})$ for $\mathbf{x} \in \partial \Omega$. $\xi_n$ is chosen to be $\delta_n/c$, where $c$ is a hyperparameter that controls the strength of the penalty. Such a choice obviously meets the assumption $\lim_{n\rightarrow\infty}\xi_n=0$ and $\lim_{n\rightarrow\infty}\frac{\delta_n^2}{\xi_n}=0$ for any $c\in \mathbb{R}$. In our experiments, we set $c=1000$ for the 1D example and $c=10$ for the 2D example, for reasons that will be explained later.
    
    The error is defined as the root mean squared error on all basis nodes:
    \[\text{Error}_{\text{1D}} = \sqrt{\frac{1}{3N+1}\sum_{j=0}^{3N}\left(u_{\text{gt}}(x_j)-u(x_j)\right)^{2}}\]
    \[\text{Error}_{\text{2D}}=\sqrt{\frac{1}{(3N+1)^2}\sum_{j_1=0}^{3N} \sum_{j_2=0}^{3N}\left(u_{\text{gt}}(x_{j_1,j_2})-u(x_{j_1,j_2})\right)^{2}}\]

    \begin{figure}[h]
        \captionsetup[subfigure]{justification=centering}
        \centering
        \subfloat[{$\Omega=[0,1], u_{\text{gt}}(x)=x^{10}$}]{\includegraphics[width=0.5\linewidth]{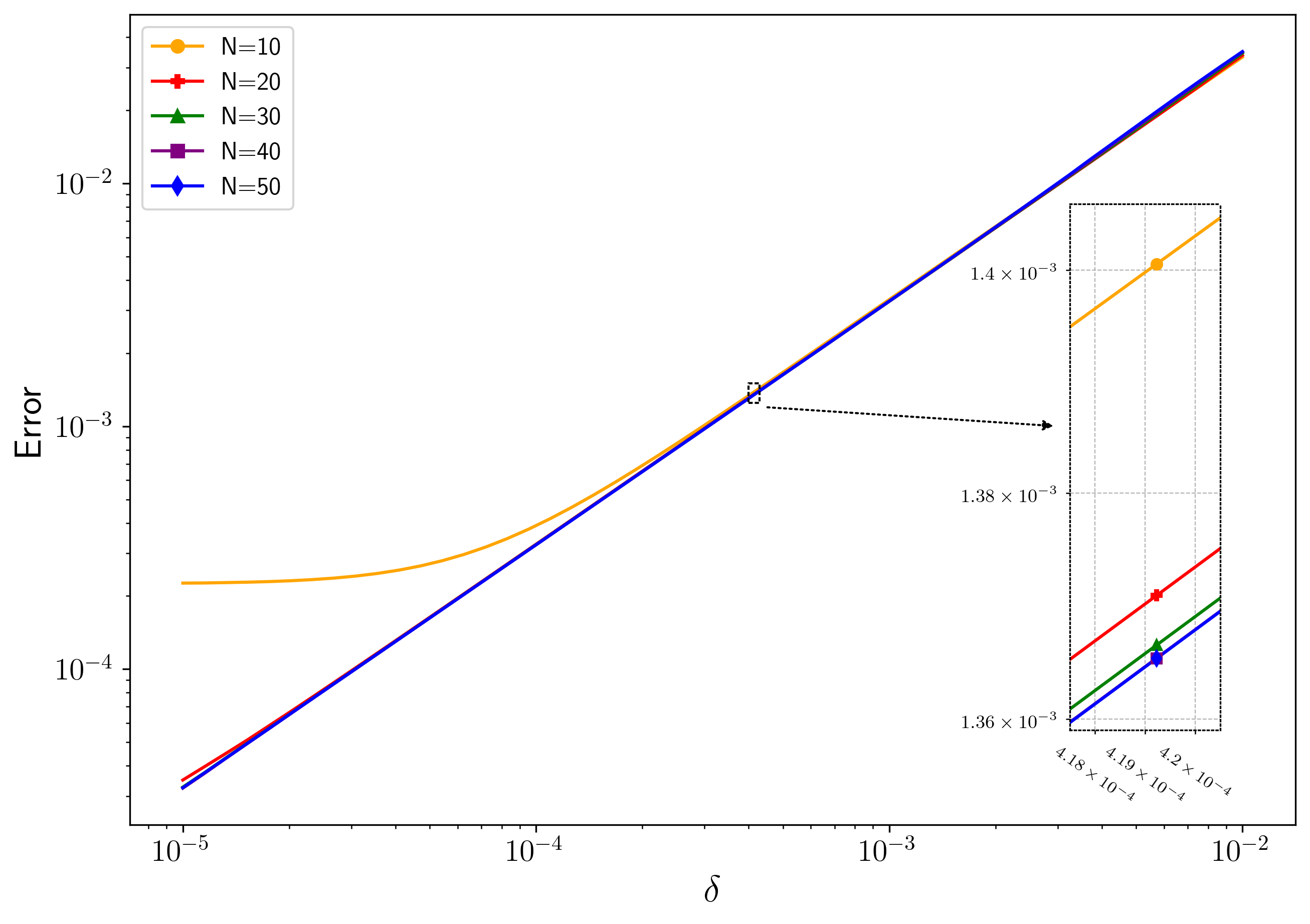}}
        \hfill
    	\subfloat[{$\Omega=[0,1]\times [0,1], u_{\text{gt}}(x_1,x_2)=x_1 \ln (x_2+1)$}]{\includegraphics[width=0.5\linewidth]{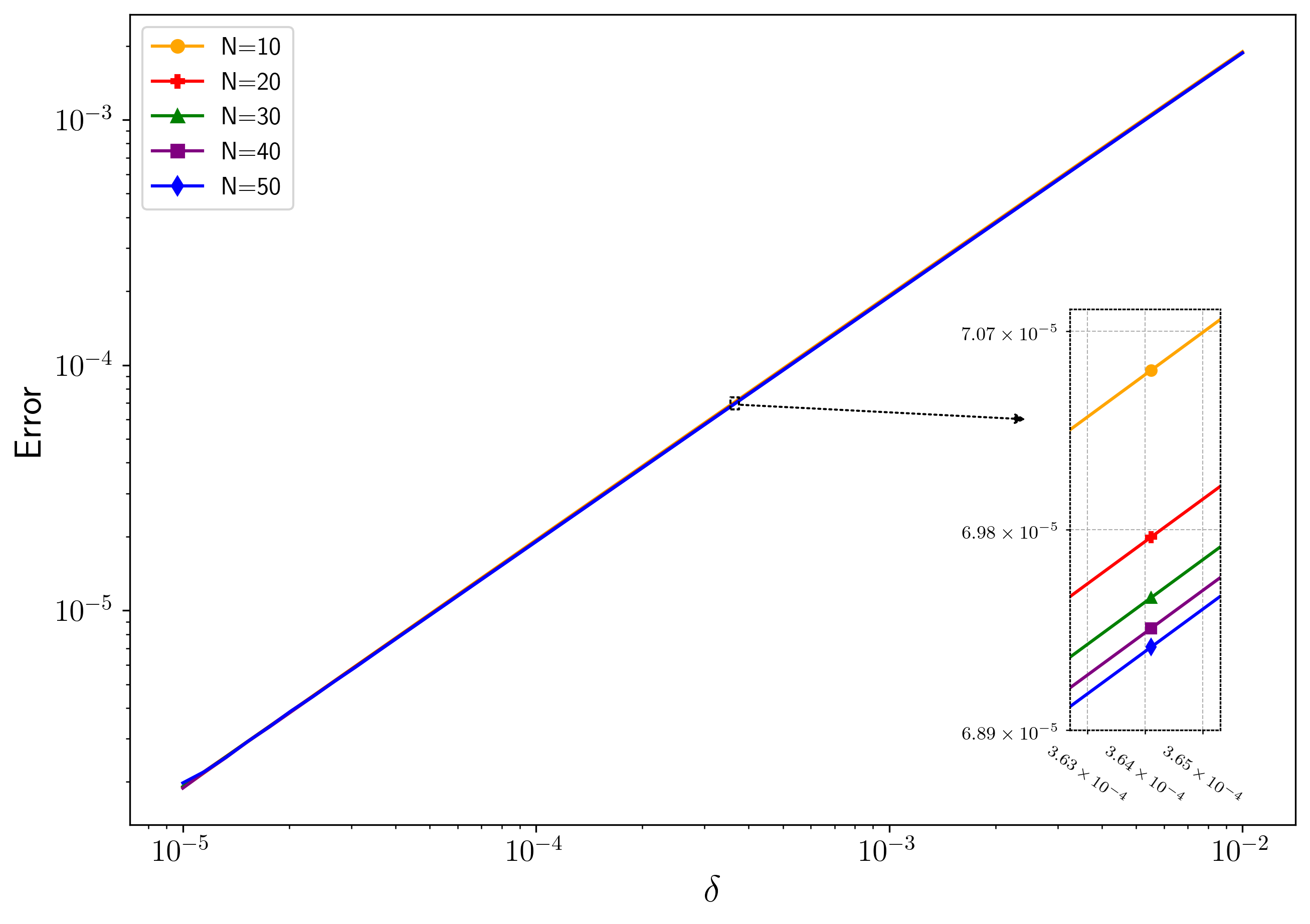}}
        \caption{Error plot with different $\delta$ in log-log scale. Lines in different colors and markers represent different $N$. A zoomed-in view of the boxed region is provided as an inset to better visualize the differences among the curves.}
        \label{fig:result1}
    \end{figure}

    The errors for different values of $\delta$ and $N$ are shown in Fig.~\ref{fig:result1}. The results validate the convergence of the nonlocal biharmonic model as $\delta$ goes to 0. In the 1D example, when $N=10$, the error is relatively large, which may be due to the insufficient accuracy of the piecewise cubic polynomial approximation for small $N$. The approximation improves as $N$ increases. { In 1D and 2D examples, the error with $N=20$ is less than $10^{-5}$ and $10^{-6}$ as shown in the inset of Fig.~\ref{fig:result1}. The reference solution is obtained with $N=50$. Notice that the error between nonlocal and local solution is around $1.37\times 10^{-3}$ and $6.9\times 10^{-5}$. This shows that $N=20$ is enough to elaborate the convergence as $\delta\rightarrow 0$. } 

    \begin{figure}[h]
        \centering
        \includegraphics[width=0.5\linewidth]{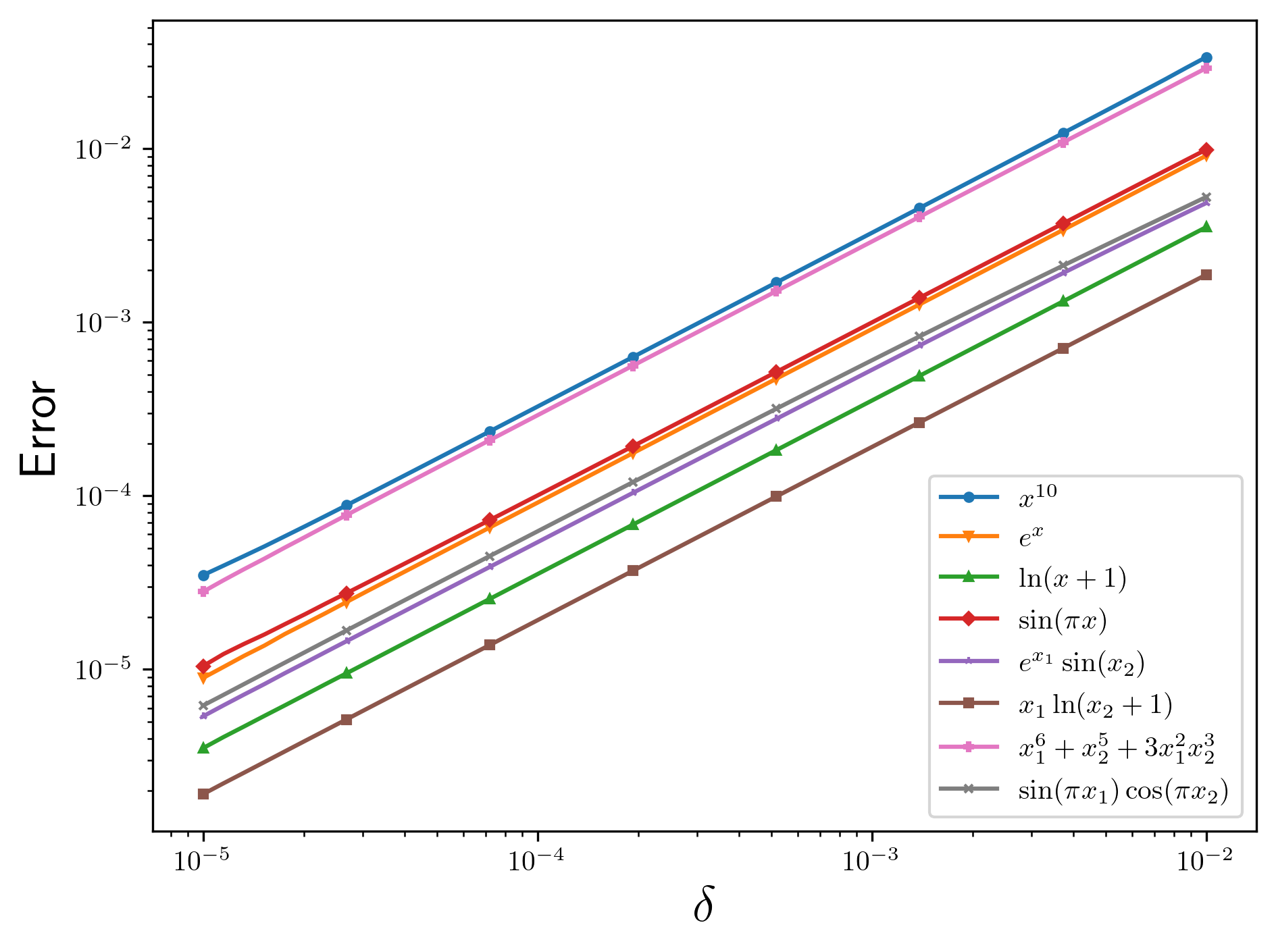}
        \caption{Error plot with $\delta$ in log-log scale. Lines in different colors and markers represent different ground-truth functions $u_{\text{gt}}$. $N=20$ is fixed. }
        \label{fig:result2}
    \end{figure}

    Moreover, the slopes of the lines in both 1D and 2D examples are very close to 1 in Fig.~\ref{fig:result1}. The convergence rate is further illustrated in Fig.~\ref{fig:result2}, where we plot the error versus $\delta$ for different examples in both 1D and 2D. The results exhibit almost same slopes very close to 1. This provides a strong evidence that the convergence rate of the proposed nonlocal biharmonic model with respect to $\delta$ is $O(\delta)$.
    
    \begin{figure}[h]
        \centering
        \subfloat[$\text{BD Error}_{\text{1D}}$]{\includegraphics[width=0.5\linewidth]{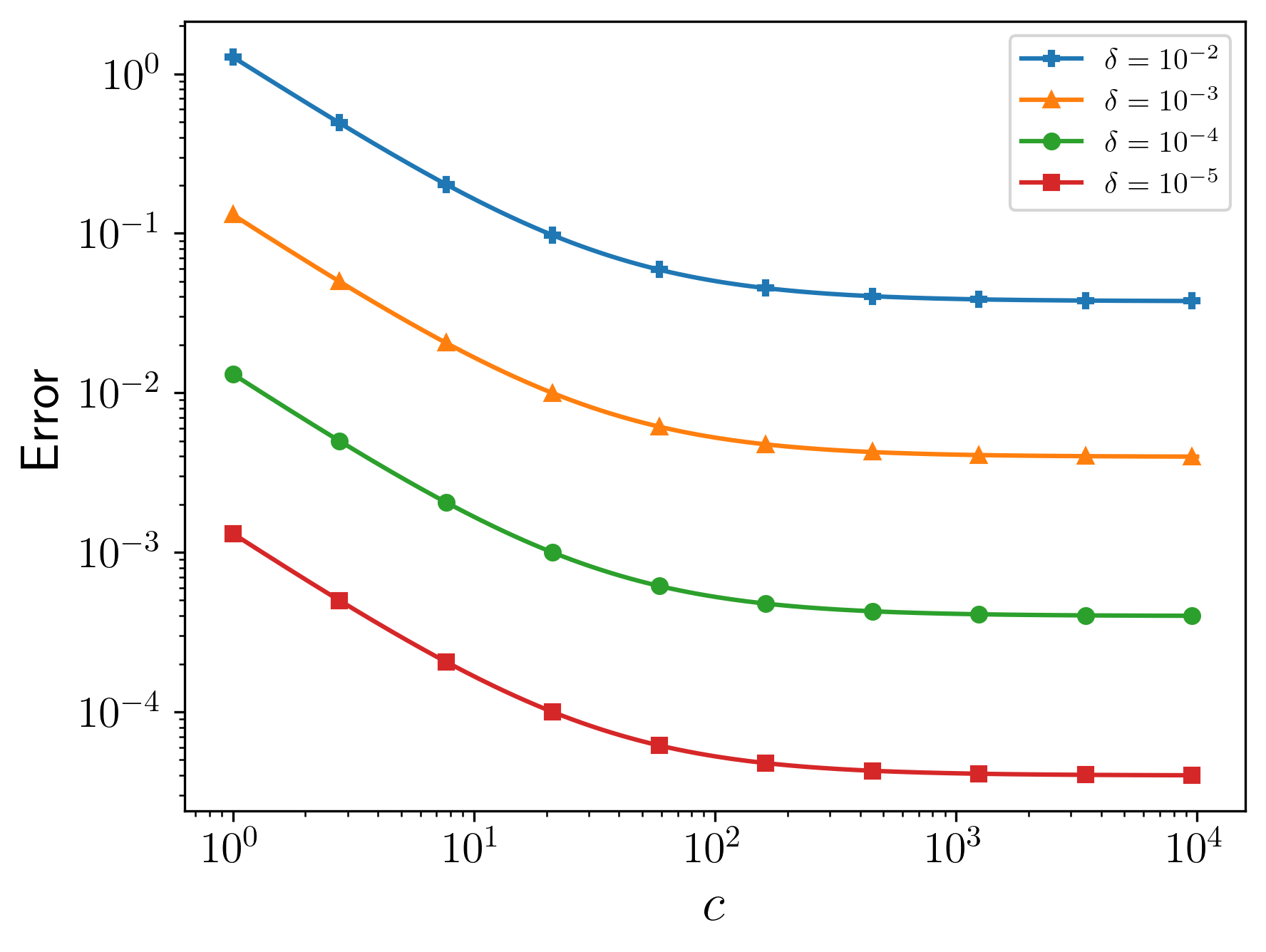}}
        \subfloat[$\text{BD $\partial_nu$ Error}_{\text{1D}}$]{\includegraphics[width=0.5\linewidth]{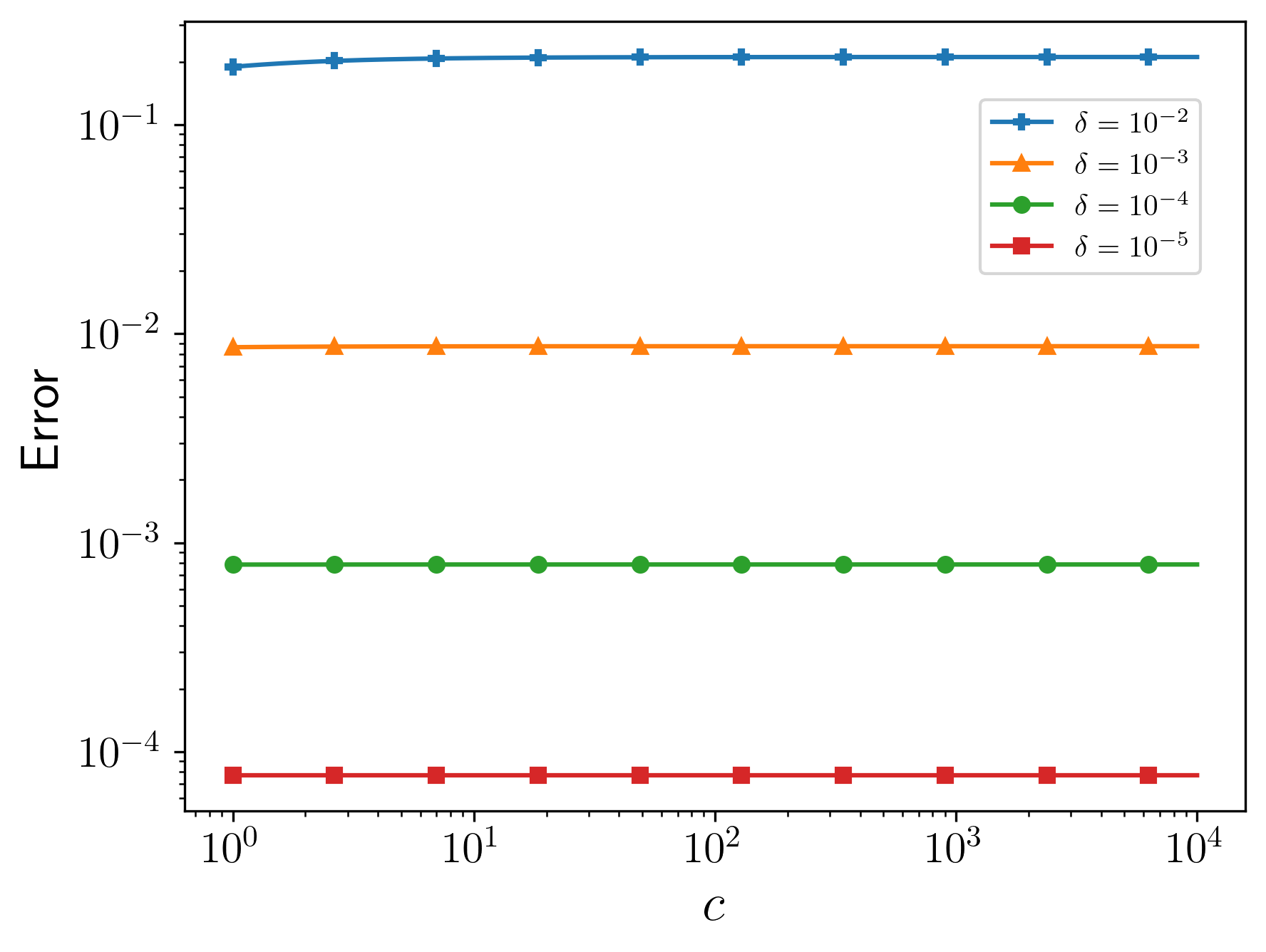}}

        \subfloat[$\text{BD Error}_{\text{2D}}$]{\includegraphics[width=0.5\linewidth]{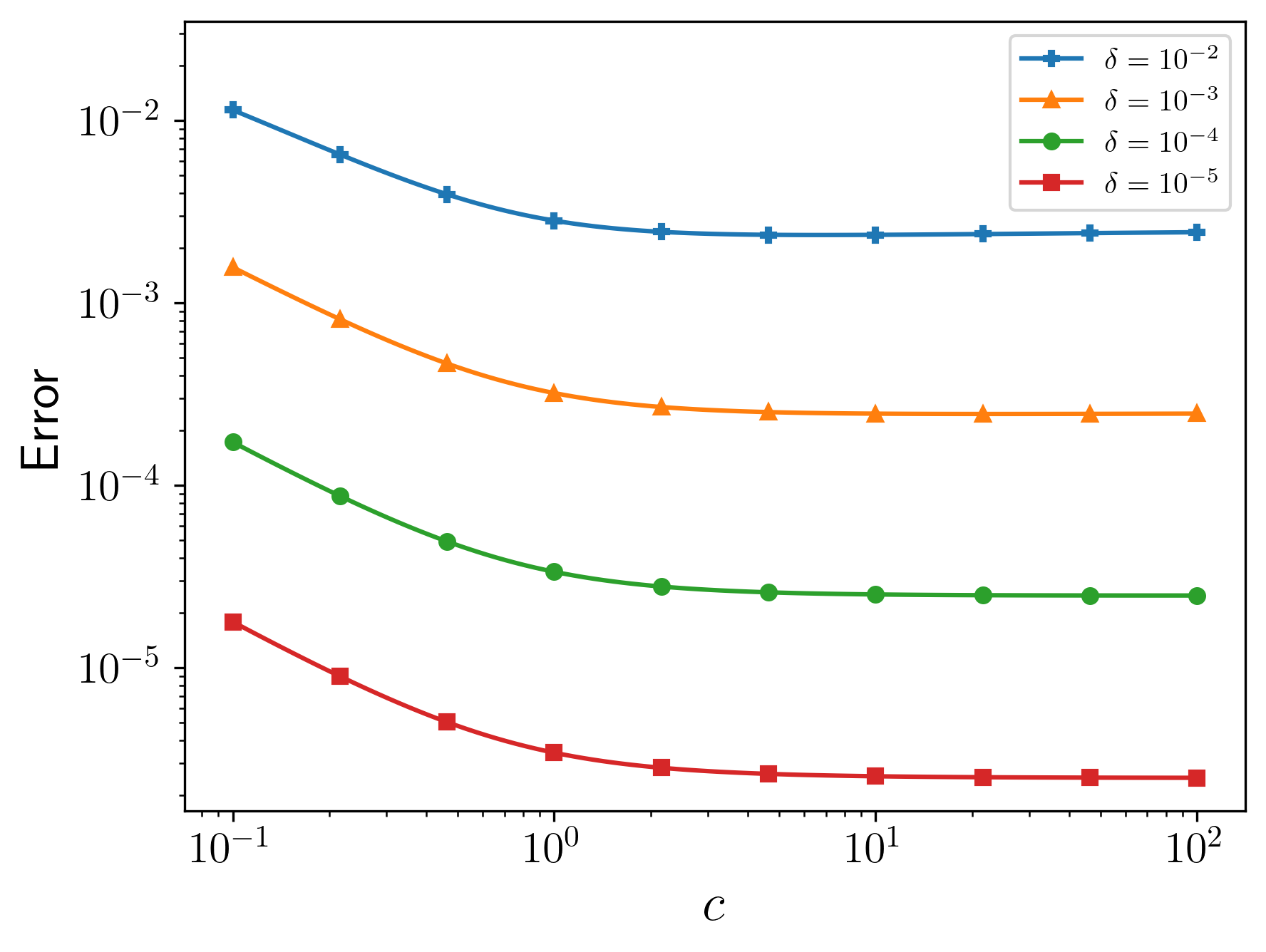}}
        \subfloat[$\text{BD $\partial_nu$ Error}_{\text{2D}}$]{\includegraphics[width=0.5\linewidth]{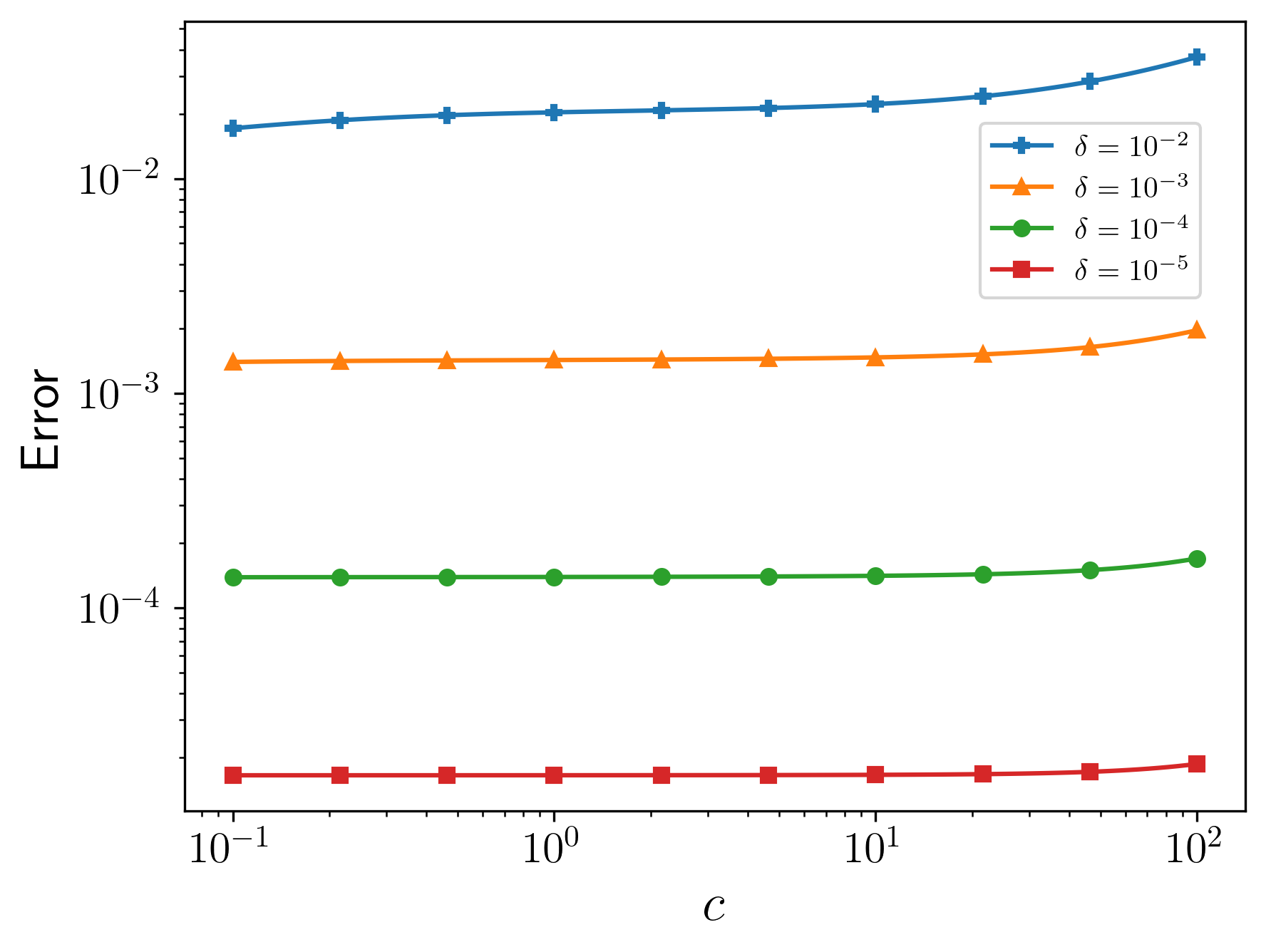}}
        \caption{Error plot with $c$ as in $\xi=\delta/c$, in log-log scale. Lines in different colors and markers represent different $\delta$. $N=20$ is fixed.}
        \label{fig:result3}
    \end{figure}

    Next, we study the effect of the parameter $\xi=\delta / c$. We fix $N = 20$ and plot the errors of both the function value and the normal derivative on the boundary with respect to $c$ for different values of $\delta$, as shown in Fig.~\ref{fig:result3}. The errors are defined as follows:
    \[\text{BD Error}_{\text{1D}} = \sqrt{\frac{1}{2}\left(\left(u_{\text{gt}}(0)-u(0)\right)^{2}+\left(u_{\text{gt}}(1)-u(1)\right)^{2} \right)}\]
    \[\text{BD $\partial_nu$ Error}_{\text{1D}} = \sqrt{\frac{1}{2}\left(\left(\frac{\partial u_{\text{gt}}}{\partial \mathbf{n}}(0)-\frac{\partial u}{\partial \mathbf{n}}(0)\right)^{2}+\left(\frac{\partial u_{\text{gt}}}{\partial \mathbf{n}}(1)-\frac{\partial u}{\partial \mathbf{n}}(1)\right)^{2} \right)}\]
    \[\text{BD Error}_{\text{2D}} = \sqrt{\frac{1}{12N} \left( \sum_{(j_1,j_2) \in \text{bd\_idx}} \left(u_{\text{gt}}(x_{j_1,j_2}) -u(x_{j_1,j_2}) \right)^{2} \right) }\]
    \[\text{BD $\partial_nu$ Error}_{\text{2D}} = \sqrt{\frac{1}{12N-4} \left( \sum_{(j_1,j_2) \in \text{bd\_idx}\backslash \text{cnr\_idx}} \left(\frac{\partial u_{\text{gt}}}{\partial \mathbf{n}}(x_{j_1,j_2}) -\frac{\partial u}{\partial \mathbf{n}}(x_{j_1,j_2}) \right)^{2} \right) }\]
    Here the boundary index set is defined as
    \[\text{bd\_idx} = \left\{ (j_1, j_2) \mid 0 \leq j_1, j_2 \leq 3N,\; j_1 = 0, 3N \text{ or } j_2 = 0, 3N \right\}\]
    and the corner index set is $\text{cnr\_idx} = \left\{ (0, 0), (0, 3N), (3N,0), (3N,3N) \right\}$. We exclude the corner points when computing the error of the normal derivative. In the one-dimensional example, we observe that as $c$ increases, that is, as the penalty becomes stronger, the error of the function value on the boundary gradually decreases and eventually stagnates when $c\ge 1000$.
    This fits the intuition very well. 
    A strong penalty tends to decrease the boundary error, but the boundary error does not decrease to zero because the error depends on $\delta$. Meanwhile, the boundary normal derivative error has little change because $\xi$ only affects the last term in $F_n(u)$, and this term does not involve the normal derivative. Thus, the effect of $\xi$ on the normal derivative error is mild. Based on these observations, we adopt $c = 1000$ in 1D experiments. 
    
    In the two-dimensional case, the phenomenon of the boundary function value error is similar to the 1D case. The only difference is that the error tends to stagnate
    when $c$ is larger than 10. For the boundary normal derivative error, it first remains nearly constant for $c\leq 10$, and begins to increase as $c$ becomes larger. Therefore, we select $c = 10$ in the 2D case. 
    
    \section{Conclusion}
    We introduce a nonlocal biharmonic model with clamped\\ boundary condition. The nonlocal model is derived based on variational approach. The well-posedness and $\Gamma$-convergence of the nonlocal model are analyzed. And the numerical experiments are conducted to validate the convergence of the nonlocal model. In the subsequent research, we will study the nonlocal biharmonic model with other boundary conditions or interface problems. 

    The numerical method also deserves more studies. The computation of the integral is simplified by utilizing that Gaussian kernel is separable in dimensions. Notice that Fourier basis or polynomial basis are also separable. For general kernels, we can approximate them using Fourier basis or polynomials. Then the integral can be computed efficiently. The methods and results will be reported in the subsequent paper. 
    
        \bibliographystyle{amsplain} 
        \bibliography{references}

\end{document}